\documentclass[12pt,a4paper,twoside]{article}

\title{\sc On the Polynomial and Exponential Decay of Eigen-Forms of Generalized Time-Harmonic Maxwell Problems}
\def\shorttitle{Polynomial and Exponential Decay of Eigen-Forms}
\def\pauthor{Dirk Pauly}

\def\mylabelonoff{off}
\def\allowdisbrk{no}

\usepackage{a4,exscale,ifthen,amsfonts,amssymb,amsmath,amscd,graphicx}
\usepackage[english]{babel}
\usepackage[mathscr]{eucal}

\author{{\sf\pauthor}}
\numberwithin{equation}{section}

\newenvironment{acknow}{{\vspace*{1cm}\noindent\bf Acknowledgements }}{}
\newcommand{\bewboxw}{\mbox{}\hfill $\square$ \\}

\newenvironment{proof}{{\noindent\bf Proof }}{\bewboxw}

\newcommand{\keywords}[1]{{\noindent\bf Key Words }#1}
\newcommand{\amsclass}[1]{{\noindent\bf AMS MSC-Classifications }#1}

\ifthenelse{\equal{\mylabelonoff}{on}}
{\newcommand{\mylabel}[1]{\label{#1}\fbox{{\rm #1}}}}{\newcommand{\mylabel}[1]{\label{#1}\makebox[0mm][]{}}}
\ifthenelse{\equal{\allowdisbrk}{yes}}
{\allowdisplaybreaks}{}

\newcommand{\paper}[7]{\bibitem{#1} #2, `#3', {\it #4}, #5, (#6), #7.}

\newcommand{\book}[6]{\bibitem{#1} #2, {\it #3}, #4, #5, (#6).}

\usepackage{amsfonts,amssymb,amsmath,amscd}
\usepackage[mathscr]{eucal}

\newcommand{\schluss}{\ifodd\value{page}\newpage\thispagestyle{empty}\makebox[0mm][]{}\color{sehrhell}.\fi\end{document}}
\newcommand{\non}{\nonumber}
\newcommand{\normabst}{\hspace{-0.4ex}}

\newcommand{\ol}[1]{\overline{#1}}

\newcommand{\ub}{\underbrace}

\newcommand{\peh}{\frac{1}{2}}

\newcommand{\meh}{{-\peh}}
\newcommand{\me}{{-1}}

\newcommand{\intom}{\int_\Omega}

\newcommand{\om}{{\Omega}}

\newcommand{\dom}{{\p\om}}

\newcommand{\ohne}{\setminus}
\newcommand{\dla}{\;d\lambda}
\newcommand{\eps}{\varepsilon}
\newcommand{\epsn}{\eps_0}
\newcommand{\epsd}{\hat{\eps}}

\newcommand{\mun}{\mu_0}
\newcommand{\mud}{\hat{\mu}}

\newcommand{\Lambdad}{\hat{\Lambda}}

\newcommand{\EH}{(E,H)}

\newcommand{\EHs}{(\tilde{E},\tilde{H})}
\newcommand{\EHD}{(E_\Delta,H_\Delta)}
\newcommand{\EHF}{(E_\calF,H_\calF)}
\newcommand{\eh}{(e,h)}

\newcommand{\FG}{(F,G)}
\newcommand{\FGd}{(\hat{F},\hat{G})}
\newcommand{\FGs}{(\tilde{F},\tilde{G})}

\newcommand{\FGD}{(F_\Delta,G_\Delta)}
\newcommand{\fg}{(f,g)}

\newcommand{\restr}[2]{\left.#1 \right|_{#2}}

\newcommand{\Abb}[5]{\begin{array}{ccccc}#1&:&#2&\longrightarrow&#3\\{}&{}&#4&\longmapsto&#5\end{array}}

\newcommand{\bds}{\begin{displaystyle}}
\newcommand{\eds}{\end{displaystyle}}
\newcommand{\beq}{\begin{equation}}
\newcommand{\eeq}{\end{equation}}
\newcommand{\beqa}{\begin{eqnarray}}
\newcommand{\eeqa}{\end{eqnarray}}
\newcommand{\beqanon}{\begin{eqnarray*}}
\newcommand{\eeqanon}{\end{eqnarray*}}
\newcommand{\bamath}{$$\begin{array}}
\newcommand{\eamath}{\end{array}$$}
\newcommand{\ba}{\begin{array}}
\newcommand{\ea}{\end{array}}
\newcommand{\bc}{\begin{center}}
\newcommand{\ec}{\end{center}}
\newcommand{\bmini}{\begin{minipage}}
\newcommand{\emini}{\end{minipage}}

\newtheorem{lem}{Lemma}[section]
\newtheorem{theo}[lem]{Theorem}
\newtheorem{cor}[lem]{Corollary}
\newtheorem{defini}[lem]{Definition}
\newtheorem{rem}[lem]{Remark}

\newcommand{\set}[2]{\{#1 \;\text{\bf :}\;  #2\}}

\newcommand{\setb}[2]{\big\{#1 \;\text{\bf :}\;  #2\big\}}

\newcommand{\rz}{\mathbb{R}}
\newcommand{\rzon}{\rz\ohne\{0\}}
\newcommand{\nz}{\mathbb{N}}
\newcommand{\zz}{\mathbb{Z}}

\newcommand{\cz}{\mathbb{C}}
\newcommand{\czp}{{\cz_+}}

\newcommand{\nzn}{{\nz_0}}

\newcommand{\rN}{{\rz^N}}

\newcommand{\rzp}{{\rz_+}}
\newcommand{\pI}{\mathbb{I}}

\newcommand{\calF}{{\mathscr{F}}}

\newcommand{\calO}{{\mathscr{O}}}

\newcommand{\calS}{{\mathscr{S}}}

\newcommand{\calM}{{\mathscr{M}}}

\newcommand{\calL}{{\mathcal{L}}}

\newcommand{\zmat}[4]{\begin{bmatrix}#1&#2\\#3&#4\end{bmatrix}}

\DeclareMathOperator{\p}{\partial}
\DeclareMathOperator{\pa}{\p^\alpha}

\DeclareMathOperator{\loes}{\calL}
\DeclareMathOperator{\loesom}{\calL_\omega}

\DeclareMathOperator{\loc}{loc}
\DeclareMathOperator{\vox}{vox}
\DeclareMathOperator{\Max}{\mathsf{Max}}

\DeclareMathOperator{\id}{Id}

\DeclareMathOperator{\ret}{Re}

\DeclareMathOperator{\supp}{supp}

\DeclareMathOperator{\rot}{rot}

\DeclareMathOperator{\pdiv}{div}

\DeclareMathOperator{\ie}{i}
\DeclareMathOperator{\e}{e}

\DeclareMathOperator{\pd}{d\hspace{-0.4ex}}
\DeclareMathOperator{\dr}{\pd r}

\newcommand{\pcoverset}[2]{\overset{#2}{\mathrm{C}}{}^{#1}}
\newcommand{\pc}[1]{\pcoverset{#1}{}}

\newcommand{\cu}{\pc{\infty}}

\newcommand{\lz}{\mathrm{L}^2}
\newcommand{\lzom}{\lz(\Omega)}
\newcommand{\Lz}[1]{\lz_{#1}}
\newcommand{\Lzom}[1]{\Lz{#1}(\Omega)}

\newcommand{\Lzt}{\Lz{t}}
\newcommand{\Lztom}{\Lz{t}(\Omega)}
\newcommand{\Lzloc}{\Lz{\loc}}
\newcommand{\Lzlocom}{\Lzloc(\Omega)}

\newcommand{\h}[3]{\overset{#3}{\mathrm{H}}{}^{#1}_{#2}}

\newcommand{\hms}{\h{m}{s}{}}

\newcommand{\hmsom}{\hms(\Omega)}

\newcommand{\qc}[3]{\overset{#3}{\mathrm{C}}{}^{#1,#2}}

\newcommand{\cqn}[1]{\qc{#1}{q}{\circ}}

\newcommand{\cqun}{\cqn{\infty}}
\newcommand{\qcom}[3]{\qc{#1}{#2}{#3}(\Omega)}

\newcommand{\cqunom}{\cqun(\Omega)}

\newcommand{\cqpe}[1]{\qc{#1}{q+1}{}}

\newcommand{\cqpeu}{\cqpe{\infty}}

\newcommand{\cqpeuom}{\cqpeu(\Omega)}

\newcommand{\qlz}[1]{\mathrm{L}^{2,#1}}
\newcommand{\qlzom}[1]{\qlz{#1}(\Omega)}
\newcommand{\lzq}{\qlz{q}}
\newcommand{\lzqom}{\qlzom{q}}
\newcommand{\qLz}[2]{\qlz{#1}_{#2}}
\newcommand{\qLzom}[2]{\qLz{#1}{#2}(\Omega)}
\newcommand{\Lzq}[1]{\qLz{q}{#1}}
\newcommand{\Lzqom}[1]{\Lzq{#1}(\Omega)}
\newcommand{\qLzs}[1]{\qLz{#1}{s}}

\newcommand{\qLzt}[1]{\qLz{#1}{t}}

\newcommand{\Lzqs}{\qLzs{q}}
\newcommand{\Lzqsom}{\Lzqs(\Omega)}
\newcommand{\Lzqt}{\qLzt{q}}
\newcommand{\Lzqtom}{\Lzqt(\Omega)}

\newcommand{\lzqpe}{\qlz{q+1}}
\newcommand{\lzqpeom}{\qlzom{q+1}}

\newcommand{\Lzqpes}{\qLzs{q+1}}
\newcommand{\Lzqpesom}{\Lzqpes(\Omega)}
\newcommand{\Lzqpet}{\qLzt{q+1}}
\newcommand{\Lzqpetom}{\Lzqpet(\Omega)}

\newcommand{\Lzqgmehom}{\Lzqom{>\meh}}

\newcommand{\qh}[4]{\overset{#4}{\mathrm{H}}{}^{#1,#2}_{#3}}

\newcommand{\qhom}[4]{\qh{#1}{#2}{#3}{#4}(\Omega)}

\newcommand{\hqmsom}{\qhom{m}{q}{s}{}}

\newcommand{\pr}[4]{{}_{#3}\overset{#4}{\mathrm{R}}{}^{#1}_{#2}}
\newcommand{\prq}[1]{\pr{q}{#1}{}{}}

\newcommand{\rqs}{\pr{q}{s}{}{}}

\newcommand{\rqt}{\pr{q}{t}{}{}}

\newcommand{\rqvox}{\pr{q}{\vox}{}{}}

\newcommand{\rqsn}{\pr{q}{s}{0}{}}
\newcommand{\rqpesn}{\pr{q+1}{s}{0}{}}

\newcommand{\rqpetn}{\pr{q+1}{t}{0}{}}

\newcommand{\ronq}[1]{\pr{q}{#1}{}{\circ}}
\newcommand{\ronqpe}[1]{\pr{q+1}{#1}{}{\circ}}
\newcommand{\ronqs}{\pr{q}{s}{}{\circ}}

\newcommand{\ronqt}{\pr{q}{t}{}{\circ}}

\newcommand{\ronqn}[1]{\pr{q}{#1}{0}{\circ}}
\newcommand{\ronqpen}[1]{\pr{q+1}{#1}{0}{\circ}}
\newcommand{\ronqsn}{\pr{q}{s}{0}{\circ}}

\newcommand{\ronqpetn}{\pr{q+1}{t}{0}{\circ}}

\newcommand{\rqom}[1]{\prq{#1}(\Omega)}

\newcommand{\rqsom}{\rqs(\Omega)}

\newcommand{\rqtom}{\rqt(\Omega)}

\newcommand{\rqvoxom}{\rqvox(\Omega)}

\newcommand{\ronqom}[1]{\ronq{#1}(\Omega)}
\newcommand{\ronqpeom}[1]{\ronqpe{#1}(\Omega)}
\newcommand{\ronqsom}{\ronqs(\Omega)}

\newcommand{\ronqtom}{\ronqt(\Omega)}

\newcommand{\ronqnom}[1]{\ronqn{#1}(\Omega)}
\newcommand{\ronqpenom}[1]{\ronqpen{#1}(\Omega)}
\newcommand{\ronqsnom}{\ronqsn(\Omega)}

\newcommand{\ronqpetnom}{\ronqpetn(\Omega)}

\newcommand{\pdi}[4]{{}_{#3}\overset{#4}{\mathrm{D}}{}^{#1}_{#2}}
\newcommand{\pdq}[1]{\pdi{q}{#1}{}{}}
\newcommand{\dqpe}[1]{\pdi{q+1}{#1}{}{}}
\newcommand{\dqs}{\pdi{q}{s}{}{}}
\newcommand{\dqpes}{\pdi{q+1}{s}{}{}}

\newcommand{\dqpet}{\pdi{q+1}{t}{}{}}

\newcommand{\dqpevox}{\pdi{q+1}{\vox}{}{}}
\newcommand{\dqn}[1]{\pdi{q}{#1}{0}{}}
\newcommand{\dqpen}[1]{\pdi{q+1}{#1}{0}{}}
\newcommand{\dqsn}{\pdi{q}{s}{0}{}}
\newcommand{\dqpesn}{\pdi{q+1}{s}{0}{}}
\newcommand{\dqtn}{\pdi{q}{t}{0}{}}

\newcommand{\dqom}[1]{\pdq{#1}(\Omega)}
\newcommand{\dqpeom}[1]{\dqpe{#1}(\Omega)}
\newcommand{\dqsom}{\dqs(\Omega)}
\newcommand{\dqpesom}{\dqpes(\Omega)}

\newcommand{\dqpetom}{\dqpet(\Omega)}

\newcommand{\dqpevoxom}{\dqpevox(\Omega)}
\newcommand{\dqnom}[1]{\dqn{#1}(\Omega)}
\newcommand{\dqpenom}[1]{\dqpen{#1}(\Omega)}

\newcommand{\dqtnom}{\dqtn(\Omega)}

\newcommand{\skp}[2]{\langle#1,#2\rangle}

\newcommand{\skpb}[2]{\big\langle#1,#2\big\rangle}

\newcommand{\norm}[1]{|\normabst|#1|\normabst|}
\newcommand{\normb}[1]{\big|\normabst\big|#1\big|\normabst\big|}

\newcommand{\qqLz}[3]{\qLz{#1,#2}{#3}}
\newcommand{\Lzqqpe}[1]{\qqLz{q}{q+1}{#1}}

\newcommand{\Lzqqpet}{\Lzqqpe{t}}
\newcommand{\Lzqqpes}{\Lzqqpe{s}}
\newcommand{\Lzqqpeloc}{\Lzqqpe{\loc}}

\newcommand{\qqLzom}[3]{\qLzom{#1,#2}{#3}}
\newcommand{\Lzqqpeom}[1]{\qqLzom{q}{q+1}{#1}}

\newcommand{\Lzqqpetom}{\Lzqqpeom{t}}
\newcommand{\Lzqqpesom}{\Lzqqpeom{s}}
\newcommand{\Lzqqpelocom}{\Lzqqpeloc(\om)}

\usepackage{algorithm,algorithmic}

\newcommand{\Fs}{\tilde{F}}
\newcommand{\Gs}{\tilde{G}}

\newcommand{\fd}{\mathbf{d}}

\renewcommand{\rot}{\pd\,}
\renewcommand{\pdiv}{\delta}
\newcommand{\Ml}{M_{\Lambda}}
\DeclareMathOperator{\tr}{tr}
\DeclareMathOperator{\sw}{sw}

\newcommand{\fs}{\tilde{f}}
\renewcommand{\fd}{\hat{f}}

\newcommand{\us}{\tilde{u}}

\newcommand{\fgs}{(\fs,\gs)}
\newcommand{\fgd}{(\fd,\gd)}

\newcommand{\FdGc}{(\Fed,\Gcd)}
\newcommand{\FcGd}{(\Fcd,\Ged)}
\newcommand{\FSGS}{(F_{\calS},G_{\calS})}
\newcommand{\EdHc}{(\Eed,\Hcd)}
\newcommand{\EcHd}{(\Ecd,\Hed)}
\renewcommand{\EH}{u}
\renewcommand{\EHD}{u_{\Delta}}
\renewcommand{\EHF}{u_{\calF}}
\newcommand{\EHell}{\EH_{\ell}}
\newcommand{\EHpiell}{\EH_{\pi\ell}}
\newcommand{\FGell}{\FG_{\ell}}
\newcommand{\FGpiell}{\FG_{\pi\ell}}
\newcommand{\EHellN}{\EHell^N}
\newcommand{\FGellN}{\FGell^N}
\newcommand{\EHellI}{\EHell^I}
\newcommand{\FGellI}{\FGell^I}
\newcommand{\EHpiellI}{\EHpiell^I}
\newcommand{\EHsell}{\EHs_{\ell}}
\newcommand{\EHspiell}{\EHs_{\pi\ell}}
\newcommand{\FGsell}{\FGs_{\ell}}
\newcommand{\FGspiell}{\FGs_{\pi\ell}}
\renewcommand{\EHs}{\us}

\renewcommand{\EdHc}{u_{\pd\;,\delta}}
\renewcommand{\EcHd}{u_{\delta,\pd}}
\renewcommand{\FG}{f}
\renewcommand{\FGD}{f_{\Delta}}
\renewcommand{\FGs}{\fs}
\renewcommand{\FGd}{\fd}
\renewcommand{\FSGS}{f_{\calS}}
\renewcommand{\FdGc}{f_{\pd\;,\delta}}
\renewcommand{\FcGd}{f_{\delta,\pd}}
\renewcommand{\eh}{u_{m}}
\renewcommand{\fg}{f_{m}}
\renewcommand{\fgd}{\fd_{m}}
\renewcommand{\fgs}{\fs_{m}}
\newcommand{\Bc}{\check{B}}
\renewcommand{\pr}[4]{{}_{#3}\overset{#4}{\mathsf{D}}{}^{#1}_{#2}}
\renewcommand{\pdi}[4]{{}_{#3}\overset{#4}{\mathsf{\Delta}}{}^{#1}_{#2}}
\renewcommand{\lz}{\mathsf{L}^2}
\renewcommand{\qlz}[1]{\mathsf{L}^{2,#1}}
\renewcommand{\qh}[4]{\overset{#4}{\mathsf{H}}{}^{#1,#2}_{#3}}
\renewcommand{\h}[3]{\overset{#3}{\mathsf{H}}{}^{#1}_{#2}}
\newcommand{\li}{\mathsf{L}^\infty}
\newcommand{\liom}{\li(\om)}
\renewcommand{\pcoverset}[2]{\overset{#2}{\mathsf{C}}{}^{#1}}
\renewcommand{\qc}[3]{\overset{#3}{\mathsf{C}}{}^{#1,#2}}
\newcommand{\vecspace}{\mathsf{V}}

\newcommand{\ronqkmeh}{\ronq{<\meh}}
\newcommand{\ronqkmehom}{\ronqkmeh(\om)}
\newcommand{\rqgmeh}{\prq{>\meh}}
\newcommand{\rqgmehom}{\rqgmeh(\om)}
\newcommand{\dqpekmeh}{\dqpe{<\meh}}
\newcommand{\dqpekmehom}{\dqpekmeh(\om)}
\newcommand{\dqpegmeh}{\dqpe{>\meh}}
\newcommand{\dqpegmehom}{\dqpegmeh(\om)}
\newcommand{\Lzqqpegmehom}{\Lzqqpeom{>\meh}}
\newcommand{\Lzqqpegehom}{\Lzqqpeom{>\peh}}

\newcommand{\omell}{\omega_{\ell}}
\newcommand{\omepill}{\omega_{\pi\ell}}

\begin{document}

\maketitle{}

\begin{abstract}
We prove polynomial and exponential decay at infinity 
of eigen-vectors of partial differential operators 
related to radiation problems for time-harmonic generalized Maxwell systems
in an exterior domain $\om\subset\rN$, $N\geq1$, 
with non-smooth inhomogeneous, anisotropic coefficients
converging near infinity with a rate $r^{-\tau}$, $\tau>1$, towards the identity.
As a canonical application we show that the corresponding eigen-values 
do not accumulate in $\rzon$
and that by means of Eidus' limiting absorption principle 
a Fredholm alternative holds true.\\
\\
\amsclass{35Q60, 78A25, 78A30}
\end{abstract}

\tableofcontents

\section{Introduction}

To establish a solution theory for time-harmonic boundary value problems in exterior domains 
it is now well known that Eidus' limiting absorption principle \cite{eidusla} is a major tool.
For this, one crucial step is to show that there are no point eigen-values,
or at least that possible point eigen-values do not accumulate and that the corresponding
eigen-spaces are finite dimensional.
The absence of non vanishing eigen-vectors can be proved by a general pattern,
which was suggested by Vogelsang \cite{volker,volkertwo} and Eidus \cite{eidussm}
and consists of the following partial results:
\begin{description}
\item[\quad step 1:] eigen-solutions decay polynomially 
\item[\quad step 2:] eigen-solutions decay exponentially
\item[\quad step 3:] eigen-solutions have compact support
\item[\quad step 4:] eigen-solutions vanish
\end{description}
These results are well known, for instance, for Helmholtz' equation including perturbations.
See \cite{agmon,eidussm,willi,kato,rellich,roze,volkertwo} and the literature cited there.
In the case of time-harmonic Maxwell's equations, steps 1 and 2 have been shown by Eidus \cite{eidussm} 
and step 4, the unique continuation property, is an old result due to Leis \cite{leistheo,leispef,leisbuch}.
The pattern was just recently completed in a sufficient manner by Bauer \cite{bauerabsence}, 
who could prove the last remaining step 3.
All these results are known for $\pc{2}$-coefficients with proper decay at infinity
except of step 1, which could have been proved even for $\li$-coefficients
by Picard, Weck and Witsch \cite{xmas}.  

In the paper at hand we address to the steps 1 and 2 for a generalized time-harmonic Maxwell problem
formulated in the language of alternating differential forms.
To show step 1 we follow closely the arguments of Picard, Weck and Witsch \cite{xmas} 
and step 2 will be proved by the methods of Eidus \cite{eidussm}.
We note that steps 3 and 4 are still open problems in our general case.
The only known result for step 4 is the case of scalar-valued $\pc{2}$-coefficients.

We consider an exterior domain, i.e. a connected open set with compact complement, 
$\om\subset\rN$, $N\geq1$, 
as a $N$-dimensional Riemannian manifold with compact boundary
and the generalized time-harmonic Maxwell equations with real frequency $\omega\neq0$
\begin{align}
\pdiv H+\ie\omega\eps E&=-\ie\eps F,\quad\rot E+\ie\omega\mu H=-\ie\mu G&&\text{in }\om,\mylabel{Maxglrotdiv}\\
\iota^* E&=0&&\text{in }\dom,\mylabel{MaxRandBedz}
\end{align}
together with the corresponding radiation condition
\begin{align}
(-1)^{qN}*\dr\wedge*H+E,\,\dr\wedge E+H\quad\text{decay at infinity.}\mylabel{MaxRadCon}
\end{align}
Here $E$, $F$ and $H$, $G$ are differential forms of rank $q$ ($q$-forms)
and $q+1$ ($(q+1)$-forms), respectively, 
and $\rot$ resp. $\pdiv=(-1)^{qN}*\rot*$ is the exterior differential resp. co-differential,
the latter acting on $(q+1)$-forms.
By $*$ we denote as usual Hodge's star operator
and by $\wedge$ the exterior product. 
$\iota:\dom\hookrightarrow\ol{\om}$ is the natural embedding of the boundary
and $\iota^*$ is the pull-back of $\iota$, i.e. the tangential trace operator.
We intend to model non-smooth, inhomogeneous and anisotropic media by linear transformations
$\eps$ and $\mu$ on $q$- and $(q+1)$-forms, respectively.

For the sake of a short notation we introduce the pars of $q$-$(q+1)$-forms
$$\EH:=(E,H),\quad\FG:=(F,G)$$
and for those the formal matrix operators
$$M:=\zmat{0}{\pdiv}{\rot}{0},\quad S:=\zmat{0}{T}{R}{0},\quad
\Lambda:=\zmat{\eps}{0}{0}{\mu},\quad\Ml:=\ie\Lambda^\me M,$$
where $R:=\dr\wedge$ and $T:=(-1)^{qN}*R*$, 
and write our problem \eqref{Maxglrotdiv}-\eqref{MaxRadCon} more compactly as
\begin{align}
(M_{\Lambda}-\omega)\EH&=\FG&&\text{in }\om,\non\\
\iota^* E&=0&&\text{in }\dom,\mylabel{shortproblem}\\
(S+1)\EH&\text{ decays at infinity}.\non
\end{align}

For the system \eqref{shortproblem} we will show polynomial and exponential decay of eigen-forms.
For the polynomial decay we can admit $\li$-coefficients $\eps$, $\mu$, 
while we need $\pc{2}$-coefficients to prove exponential decay.
In both cases the coefficients must converge at infinity with a rate $r^{-\tau}$, $\tau>1$, 
towards homogeneous and isotropic coefficients.

The main tool to handle irregular coefficients is a decomposition lemma, 
which allows us to prove the polynomial decay of eigen-forms
by reduction to the similar result known for the scalar Helmholtz equation. 
The keys to this decomposition lemma are weighted Hodge-Helmholtz decompositions, 
i.e. decompositions into irrotational and solenoidal forms, in the whole space case
and a well known procedure to decouple the electric and magnetic form 
by discussing a second order elliptic system. 
To illustrate this calculation let us look at \eqref{shortproblem}
in the homogeneous case. Applying $M_{\id}+\omega$ yields
\beq(M^2+\omega^2)\EH=0.\mylabel{calcsecorderz}\eeq
By \eqref{shortproblem} $E$ is solenoidal, i.e. $\pdiv E=0$, and $G$ irrotational, i.e. $\rot H=0$, 
since $\pdiv\pdiv=0$, $\rot\rot=0$. From $\Delta=\rot\pdiv+\pdiv\rot$, 
where the Laplacian acts on each Euclidian component, 
we get the identity $M^2\EH=(\pdiv\rot E,\rot\pdiv H)=\Delta\EH$ and
finally \eqref{calcsecorderz} turns to the componentwise Helmholtz equation
\beq(\Delta+\omega^2)\EH=0.\mylabel{calcsecorderd}\eeq

The polynomial decay of eigen-forms together with an a priori estimate 
for the solutions corresponding to non-real frequencies
is sufficient to prove a Fredholm alternative for \eqref{shortproblem}
utilizing the limiting absorption principle invented by Eidus \cite{eidusla}.
Moreover, we get at most finite dimensional eigen-spaces for possible eigen-values
but these can not accumulate in $\rzon$. 

\section{Definitions and preliminaries}

For later purpose let us fix $r_{0}>0$, such that $\rN\ohne\om\subset B_{r_{0}}$,
where $B_{\theta}$ denotes the open ball of radius $\theta$ centered at the origin.
We also define the exterior of the closed ball $\Bc_{\theta}:=\rN\ohne\ol{B_{\theta}}$
and the sphere $S_{\theta}$, both of radius $\theta$.  

Using the weight function
$$\rho:=(1+r^2)^{1/2},\quad r(x):=x,$$
we introduce for $m\in\nzn$ and $s\in\rz$ the weighted scalar Sobolev spaces
$$\hmsom:=\setb{\psi\in\Lzlocom}{\rho^s\p^\alpha\psi\in\lzom\quad\forall\,|\alpha|\leq m}.$$
In $\om$ we have a canonical global chart, the identity, and thus, 
$\om$ becomes naturally a $N$-dimensional smooth Riemannian manifold 
with Cartesian coordinates $\{x_1,\dots,x_N\}$. 
For alternating differential forms of rank $q\in\zz$ ($q$-forms)
we define componentwise partial derivatives $\p^\alpha\Phi=(\p^\alpha\Phi_I)\pd x^I$, 
if $\Phi=\Phi_I\pd x^I$ (sum convention!), where $I$ are ordered multi-indices of length $q$.
Then, for $m\in\nzn$ and $s\in\rz$ we define weighted Sobolev spaces $\hqmsom$ of $q$-forms as well.
Equipped with their natural scalar products all these spaces become Hilbert spaces.
For $m=0$ we also utilize the notation $\Lzqsom:=\qhom{0}{q}{s}{}$ 
and for pairs of forms we introduce product spaces like
$$\Lzqqpesom:=\Lzqsom\times\Lzqpesom.$$  
In the special case $s=0$ we neglect the index $0$
and we have in $\lzqom=\qhom{0}{q}{}{}=\qhom{0}{q}{0}{}$ the scalar product
$$\skp{\Phi}{\Psi}_{\Lzqom{}}=\intom\Phi\wedge*\ol{\Psi}
=\intom*\skp{\Phi}{\Psi}_q=\intom\skp{\Phi}{\Psi}_q\dla=\intom\Phi_I\ol{\Psi}_I\dla$$
for $\Phi=\Phi_I\pd x^I$, $\Psi=\Psi_I\pd x^I\in\Lzqom{}$.
Here $\lambda$ denotes Lebesgue's measure and $\skp{\,\cdot\,}{\,\cdot\,}_q$ the pointwise scalar product on $q$-forms.

By Stokes' theorem and the product rule the exerior derivative and co-derivative 
are formally skew-adjoint to each other, i.e. 
$$\skp{\rot\Phi}{\Psi}_{\lzqpeom}=-\skp{\Phi}{\pdiv\Psi}_{\lzqom}
\quad\forall\,(\Phi,\Psi)\in\qcom{\infty}{q,q+1}{\circ},$$
which gives rise to weak definitions of $\rot$ and $\pdiv$. 
Here we denote the vector space of all smooth $q$-forms with compact support in $\om$
by $\cqunom$. We note that still $\rot\rot=0$,
$\pdiv\pdiv=0$ and $\rot\pdiv+\pdiv\rot=\Delta$ hold true in the weak sense.
Furthermore, for $s\in\rz$ we introduce some special weighted Sobolev spaces suited for Maxwell's equations
\begin{align*}
\rqsom&:=\setb{\Phi\in\Lzqsom}{\rot\Phi\in\Lzqpesom},\\
\dqsom&:=\setb{\Phi\in\Lzqsom}{\pdiv\Phi\in\qLzom{q-1}{s}}.
\end{align*}
Equipped with their natural graph norms these are Hilbert spaces as well. 
To generalize the homogeneous boundary condition we introduce 
$\ronqsom$ as the closure of $\cqun(\Omega)$ in the norm of $\rqsom$.
Utilizing Stokes' theorem we see that in fact the homogeneous boundary condition $\iota^*E=0$ 
is generalized in $\ronqsom$. 
The spaces $\rqsom$, $\dqsom$ and even $\ronqsom$ are invariant under multiplication
with bounded smooth functions $\varphi$.
A subscript $0$ at the lower left corner indicates vanishing exterior derivative resp. co-derivative,
for instance, $\ronqsnom=\setb{\Phi\in\ronqsom}{\rot\Phi=0}$.
The indices $\loc$ and $\vox$ refer as usual to local integrability and compact supports, respectively.
If the whole space $\om=\rN$ is under consideration, we omit the dependence on the domain
and write simply, for example, $\rqsn:=\rqsn(\rN)$. 
Moreover, for weighted Sobolev spaces $\vecspace_t$, $t\in\rz$, we define
$$\vecspace_{<s}:=\bigcap_{t<s}\vecspace_t,\quad\vecspace_{>s}:=\bigcup_{t>s}\vecspace_t,\quad s\in\rz.$$

Now let us introduce the properties of our transformations $\eps$, $\mu$, $\Lambda$:

\begin{defini}\mylabel{transdefi}
Let $\tau\geq0$. We call a transformation $\eps$ $\tau$-{\sf admissible}, if
\begin{itemize}
\item[\bf(i)] $\eps(x)$ is a linear transformation on $q$-forms for all $x\in\om$;
\item[\bf(ii)] $\eps$ possesses $\liom$-coefficients, i.e. the matrix representation of
$\eps$ corresponding to the canonical basis (and then for every chart basis)
has $\liom$-entries;
\item[\bf(iii)] $\eps$ is symmetric, i.e.
$$\forall\,\Phi,\Psi\in\lzqom\quad\skp{\eps\Phi}{\Psi}_{\lzqom}=\skp{\Phi}{\eps\Psi}_{\lzqom};$$
\item[\bf(iv)] $\eps$ is uniformly positive definite, i.e.
$$\exists\,c>0\,\forall\,\Phi\in\lzqom\quad\skp{\eps\Phi}{\Phi}_{\lzqom}\geq c\norm{\Phi}_{\lzqom}^2;$$
\item[\bf(v)] $\eps$ is asymptotically the identity, i.e.
$\eps=\epsn\id+\hat{\eps}$ with $\epsn\in\rzp$ and $\hat{\eps}=\calO(r^{-\tau})$ as $r\to\infty$.
\end{itemize}
Moreover, for $n\in\nzn$ we call $\eps$ $\tau$-$\pc{n}$-{\sf admissible}, if $\eps$ is $\tau$-admissible and
\begin{itemize}
\item[\bf(vi)] $\epsd\in\pc{n}(\Bc_{r_0})$ with bounded derivatives, 
which means that the matrix representation of $\epsd$ corresponding to the canonical basis
(and then for every chart basis) has $\pc{n}(\Bc_{r_0})$-entries and all derivatives are bounded.
\end{itemize}
We call $\tau$ the {\sf order of decay} of the perturbation $\hat{\eps}$.
\end{defini}

We remark that by a transformation $\tilde{x}:=\alpha x$, $\tilde{H}:=\beta H$ 
we may assume with loss of generality $\epsn=\mun=1$ throughout this paper.

Finally, we note that the multiplication operators $R$, $T$ and $S$
are related to the differential operators $\rot$, $\pdiv$, $M$ through the following formulas:
$$C_{D,\varphi(r)}=\varphi'(r)X,\quad
\calF D=\ie rX\calF,\quad D\calF=-\ie\calF rX$$ 
Here $(D,X)\in\big\{(\rot,R),(\pdiv,T),(M,S)\big\}$
and $C_{A,B}$ denotes the commutator of two operators $A$, $B$, i.e.
$C_{D,\varphi(r)}=D\varphi(r)-\varphi(r)D$,
where $\varphi$ is a smooth function on $\rz$. 
Furthermore, $\calF$ denotes the Fourier transformation on $q$-forms in $\rN$
(componentwise in Euclidean coordinates), which is an unitary mapping on $\lzq$.

\section{A decomposition lemma}

The following decomposition lemma is essential
and allows us to transfer results known from Helmholtz' equation to Maxwell's equations
without any further regularity assumptions.
To use results from Weck and Witsch \cite{sphharm} we set
$$\pI:=\set{n+N/2,1-n-N/2}{n\in\nzn}.$$

\begin{lem}\mylabel{decomplemma}
Let $\Lambda$ be $\tau$-admissible with order of decay $\tau\geq0$.
Furthermore, let $K$ be a compact subset of $\cz\ohne\{0\}$, $\omega\in K$, $t,s\in\rz$ 
with $0\leq s\in\rz\ohne\pI$ and $t\leq s\leq t+\tau$.
Let $\theta\geq r_0$ and $\varphi:=\eta(r/\theta)$, 
where $\eta\in\cu(\rz)$ is a cut-off function supported in $[1,\infty)$ 
and constantly equal to $1$ in $[2,\infty)$. 
Moreover, let $\EH\in\rqtom\times\dqpetom$ satisfy
$$(\Ml-\omega)\EH=:\FG\in\Lzqqpesom.$$
Then,
$$\FGd:=-\ie\varphi\Lambda\FG+(C_{M,\varphi}-\ie\omega\varphi\Lambdad)\EH\in\Lzqqpes$$
and by decomposing
$$\FGd=:\FdGc+\FcGd+\FSGS\in
(\rqsn\times\dqpesn)\dot{+}(\dqsn\times\rqpesn)\dot{+}\calS^{q,q+1}_s$$
according to \cite[Theorem 4]{sphharm}
$$\FGs:=\FcGd+\frac{\ie}{\omega}M\FSGS\in\dqsn\times\rqpesn$$
holds. Then $\EH$ can be decomposed into
$$\EH=(1-\varphi)\EH+\EdHc+\EHF+\EHD$$
and there exist generic constants $c>0$, which are independent of $\EH$, $\FG$ or $\omega$, such that
\begin{itemize}
\item[\rm\bf (i)] $\bds(1-\varphi)\EH\in\rqvoxom\times\dqpevoxom\eds$ and for all $\tilde{t}\in\rz$
$$\normb{(1-\varphi)\EH}_{\rqom{\tilde{t}}\times\dqpeom{\tilde{t}}}
\leq c\big(\norm{\FG}_{\Lzqqpesom}+\norm{\EH}_{\Lzqqpeom{s-\tau}}\big);$$
\item[\rm\bf (ii)] $\bds\EdHc:=-\frac{\ie}{\omega}(\FdGc+\FSGS)\in\rqs\times\dqpes\eds$ and
$$\norm{\EdHc}_{\rqs\times\dqpes}\leq c\norm{\FGd}_{\Lzqqpes};$$
\item[\rm\bf (iii)] $\EHF:=\calF^\me\big((1+r^2)^\me(1-\ie rS)\calF\FGs\big)
\in\qh{1}{q,q+1}{s}{}\cap(\dqsn\times\rqpesn)$ and
$$\norm{\EHF}_{\qh{1}{q,q+1}{s}{}}\leq c\norm{\FGs}_{\Lzqqpes};$$
\item[\rm\bf (iv)] $\EHD:=\EcHd-\EHF\in\qh{2}{q,q+1}{t}{}\cap(\dqtn\times\rqpetn)$ and
$$\norm{\EHD}_{\qh{2}{q,q+1}{\tilde{t}}{}}\leq 
c\big(\norm{\EHD}_{\Lzqqpe{\tilde{t}}}+\norm{\EHF}_{\qh{1}{q,q+1}{\tilde{t}}{}}\big)$$
for all $\tilde{t}\leq t$, where
$$\EcHd:=-\frac{\ie}{\omega}(\FcGd-M\varphi\EH)\in\qh{1}{q,q+1}{t}{}\cap(\dqtn\times\rqpetn).$$
\end{itemize}
These forms solve 
\begin{align*}
(M+\ie\omega)\varphi\EH&=\FGd,\quad
(M+\ie\omega)\EcHd=\FGs,\quad
(M+\ie\omega)\EHD=(1-\ie\omega)\EHF,\\
(M+1)\EHF&=\FGs
\intertext{and}
(\Delta+\omega^2)\EHD&=-(1+\omega^2)\EHF+(1-\ie\omega)\FGs.
\end{align*}
Moreover, the estimates 
\begin{align*}
\norm{\FGs}_{\Lzqqpes}
&\leq c\norm{\FGd}_{\Lzqqpes},\\
\norm{\FGd}_{\Lzqqpes}
&\leq c\big(\norm{\FG}_{\Lzqqpesom}
+\norm{\EH}_{\Lzqqpeom{s-\tau}}\big),\\
\norm{\EH}_{\rqom{\tilde{t}}\times\dqpeom{\tilde{t}}}
&\leq c\big(\norm{\FG}_{\Lzqqpesom}
+\norm{\EH}_{\Lzqqpeom{s-\tau}}+\norm{\EHD}_{\Lzqqpe{\tilde{t}}}\big)
\intertext{and}
\normb{(\Delta+\omega^2)\EHD}_{\Lzqqpes}
&\leq c\big(\norm{\FG}_{\Lzqqpesom}
+\norm{\EH}_{\Lzqqpeom{s-\tau}}\big)
\intertext{as well as}
\normb{(M-\ie\lambda S)\EH}_{\Lzqqpeom{\tilde{t}}}
&\leq c\Big(\norm{\FG}_{\Lzqqpesom}+\norm{\EH}_{\Lzqqpeom{s-\tau}}
+\normb{(M-\ie\lambda S)\EHD}_{\Lzqqpe{\tilde{t}}}\Big)
\end{align*}
hold for all $\tilde{t}\leq t$ and uniformly in $\lambda\in K$, $\EH$ and $\FG$.
\end{lem}

\begin{proof}
Obviously, $(\Ml-\omega)\EH=\FG\in\Lzqqpesom$ is equivalent to
$$(M+\ie\omega\Lambda)\EH=-\ie\Lambda\FG\in\Lzqqpesom.$$
Since $\varphi\EH\in\rqt\times\dqpet$ we have 
$$M\varphi\EH=\varphi M\EH+C_{M,\varphi}\EH
=-\ie\omega\varphi\Lambda\EH-\ie\varphi\Lambda\FG+C_{M,\varphi}\EH$$
and thus 
\beq(M+\ie\omega)\varphi\EH=\FGd\in\Lzqqpes,\mylabel{gleins}\eeq
since $C_{M,\varphi}=\eta'(r/\theta)\theta^\me S$ is compactly supported and $t+\tau\geq s$.
We rewrite \eqref{gleins} in the form
$$\ie\omega\varphi\EH=-M\varphi\EH+\FcGd+\FdGc+\FSGS$$
and note
\begin{align*}
\EcHd&=-\frac{\ie}{\omega}(\FcGd-M\varphi\EH)
\in(\rqt\cap\dqtn)\times(\dqpet\cap\rqpetn)\subset\qh{1}{q,q+1}{t}{},\\
\EdHc&=-\frac{\ie}{\omega}(\FdGc+\FSGS)
\in\rqs\times\dqpes
\end{align*}
with $\varphi\EH=\EcHd+\EdHc$ and by regularity, 
e.g. \cite[Lemma 4.2(i)]{kuhnpaulyreg}.
(For $s<N/2$ we even have $\FSGS=0$.) 
Moreover, $\EcHd$ solves
$$M\EcHd=M\varphi\EH-M\EdHc=-\ie\omega\EcHd+\FcGd+\frac{\ie}{\omega}M\FSGS,$$
i.e. $(M+\ie\omega)\EcHd=\FGs\in\dqsn\times\rqpesn$.
Now, to define $(M+1)^\me\FGs$ by the Fourier transformation we put
$$\EHF=\calF^\me\big((1+r^2)^\me(1-\ie rS)\calF\FGs\big).$$
Then, $\EHF\in\Lzqqpe{}$ as well as $\calF\EHF\in\Lzqqpe{1}$ 
are implied by $\calF\FGs\in\Lzqqpe{}$. Hence,
$$\EHF\in\qh{1}{q,q+1}{}{}.$$
From $\FGs\in\Lzqqpes$ we get by definition $\calF\FGs\in\qh{s}{q,q+1}{}{}$.
The components of $\calF\EHF$ arise from those of $\calF\FGs$ 
by multiplication with bounded $\cu$-functions. Thus, also
$$\calF\EHF\in\qh{s}{q,q+1}{}{}$$
follows; see e.g. Wloka \cite[p. 71, Lemma 3.2]{wloka}.
Again, by definition, $\EHF\in\Lzqqpes$ and we obtain the estimate
$$\norm{\EHF}_{\Lzqqpes}\leq c\norm{\FGs}_{\Lzqqpes}.$$
Since $M\calF^\me=\ie\calF^\me rS$ we compute
$$(1+r^2)\calF(M+1)\EHF=(1+\ie rS)(1-\ie rS)\calF\FGs=(1+r^2S^2)\calF\FGs.$$
$\pdiv\tilde{F}=0$ and $\rot\tilde{G}=0$ 
imply $T\calF\tilde{F}=0$ and $R\calF\tilde{G}=0$, respectively,
and therefore, using $RT+TR=1$
$$S^2\calF\FGs=\zmat{TR}{0}{0}{RT}\calF\FGs=\calF\FGs$$
holds. Hence,
$$\calF(M+1)\EHF=\calF\FGs\text{, i.e.}\quad(M+1)\EHF=\FGs.$$
Besides, we have $\EHF\in(\rqs\cap\dqsn)\times(\dqpes\cap\rqpesn)$ and thus,
$$\EHF\in\qh{1}{q,q+1}{s}{}\cap(\dqsn\times\rqpesn)$$
again by regularity.
Considering
$$\EHD=\EcHd-\EHF\in\qh{1}{q,q+1}{t}{}\cap(\dqtn\times\rqpetn)$$
we calculate
$$(M+\ie\omega)\EHD=(1-\ie\omega)\EHF.$$
Once more by regularity we even obtain $\EHD\in\qh{2}{q,q+1}{t}{}\cap(\dqtn\times\rqpetn)$
and we compute
\begin{align*}
(\Delta+\omega^2)\EHD&=(M-\ie\omega)(M+\ie\omega)\EHD
=(1-\ie\omega)(M-\ie\omega)\EHF\\
&=-(1+\omega^2)\EHF+(1-\ie\omega)\FGs.
\end{align*}
Finally, we achieve the asserted estimates from the regularity result and the continuity
of the projections in $\Lzqs$ onto $\rqsn$, $\dqsn$ resp. $\calS^q_s$ mentioning that
$$\norm{M\FSGS}\leq c\norm{\FSGS}$$
holds in any norm since $\calS^{q,q+1}_s$ is finite dimensional and $M$ linear.
\end{proof}

\section{Polynomial decay}

First, we need a trivial but useful technical lemma.

\begin{lem}\mylabel{normtrick}
For all $t,\tilde{t}\in\rz$ with $\tilde{t}<t$ and all $\vartheta>0$ 
there exist constants $c,\theta>0$, such that
$$\norm{\psi}_{\Lzom{\tilde{t}}}
\leq c\norm{\psi}_{\lz(\Omega\cap B_{\theta})}+\vartheta\norm{\psi}_{\Lztom}$$
holds for all $\psi\in\Lztom$.
\end{lem}

\begin{proof}
For sufficient large $\theta>0$ we get from $\tilde{t}-t<0$
$$\norm{\psi}_{\Lzom{\tilde{t}}}^2
=\norm{\rho^{\tilde{t}}\psi}_{\lz(\Omega\cap B_{\theta})}^2
+\norm{\rho^{\tilde{t}-t}\psi}_{\Lzt(\Bc_{\theta})}^2
\leq c_{\om,\tilde{t},\theta}\norm{\psi}_{\lz(\Omega\cap B_{\theta})}^2
+(1+\theta^2)^{\tilde{t}-t}\norm{\psi}_{\Lzt(\Bc_{\theta})}^2.$$
Thus, $\bds\lim_{\theta\to\infty}(1+\theta^2)^{\tilde{t}-t}=0\eds$ completes the proof.
\end{proof}

Our decomposition lemma implies:

\begin{theo}\mylabel{polyabklmax}
Let $\Lambda$ be $\tau$-admissible with $\tau>1$.
Moreover, let $I\subset\rz_{\pm}$ be a closed interval and 
$\omega\in I$ as well as $1/2<s\in\rz\ohne\pI$. If
$$\EH\in\rqgmehom\times\dqpegmehom$$
is a solution of Maxwell's equation
$$(\Ml-\omega)\EH=:\FG\in\Lzqqpesom,$$
then $\EH\in\rqom{s-1}\times\dqpeom{s-1}$ and there exist constants $c,\theta>0$ independent of $\EH$, $\FG$ or $\omega$, such that
$$\norm{\EH}_{\rqom{s-1}\times\dqpeom{s-1}}
\leq c\big(\norm{\FG}_{\Lzqqpesom}
+\norm{\EH}_{\Lzqqpe{}(\Omega\cap B_{\theta})}\big).$$
\end{theo}

\begin{proof}
Let $t>-1/2$ and $\EH\in\rqtom\times\dqpetom$ with $t<s-1$. 
Without loss of generality we may assume $t+1<s<t+\tau$.
Otherwise, we replace $t$ and $s$ by $t_k:=t+k\alpha$ and $s_k:=t+1+(k+1)\alpha\leq s$, $k=0,\dots$,
with $\alpha:=(\tau-1)/2>0$ and obtain the assertions after finitely many $\alpha$-steps.

Decomposing $\EH$ by Lemma \ref{decomplemma} we get solutions
$\EHD\in\qh{2}{q,q+1}{t}{}$ of Helmholtz' equation in $\rN$
$$(\Delta+\omega^2)\EHD=:\FGD\in\dqsn\times\rqpesn.$$
A componentwise application of \cite[Lemma 5]{linelae}
yields $\EHD\in\qh{2}{q,q+1}{s-1}{}$ and with a
constant $c>0$ independent of $\EHD$, $\FGD$ or $\omega$ we have
$$\norm{\EHD}_{\qh{2}{q,q+1}{s-1}{}}
\leq c\big(\norm{\FGD}_{\Lzqqpes}
+\norm{\EHD}_{\Lzqqpe{s-2}}\big).$$
Furthermore, from Lemma \ref{decomplemma} we have $\EH\in\dqom{s-1}\times\dqpeom{s-1}$ and the estimate 
\begin{align*}
\norm{\EH}_{\rqom{s-1}\times\dqpeom{s-1}}
&\leq c\big(\norm{\EHD}_{\Lzqqpe{s-1}}
+\norm{\FG}_{\Lzqqpesom}
+\norm{\EH}_{\Lzqqpeom{s-\tau}}\big)\\
&\leq c\big(\norm{\FGD}_{\Lzqqpes}
+\norm{\FG}_{\Lzqqpesom}
+\norm{\EH}_{\Lzqqpeom{s-\tau}}\big)\\
&\leq c\big(\norm{\FG}_{\Lzqqpesom}
+\norm{\EH}_{\Lzqqpeom{s-\tau}}\big),
\end{align*}
where we assumed without loss of generality $\tau<2$.
Since $s-\tau<s-1$ the assertion follows now by Lemma \ref{normtrick}.
\end{proof}

\begin{rem}\mylabel{bempolyabklmax}
Let the assumptions of Theorem \ref{polyabklmax} be satisfied.
If $\FG\in\Lzqqpesom$ for all $s\in\rz$, we get
$$\EH\in\rqsom\times\dqpesom$$
for all $s\in\rz$. 
This holds, for instance, if $\FG$ is exponentially decaying or even compactly supported.
\end{rem}

\section{Exponential decay}

Using the 'partial integration' technique introduced by Eidus \cite{eidussm}
for the classical Maxwell equations we will prove:

\begin{theo}\mylabel{expabkl}
Let $\omega\in\rzon$ and $\Lambda$ be $\tau$-$\pc{2}$-admissible
with $\tau>1$. If
$$\EH\in\rqgmehom\times\dqpegmehom$$
is a solution of Maxwell's equation
$$(\Ml-\omega)\EH=:\FG\in\e^{tr}\qh{2}{q,q+1}{}{}(\Bc_{r_{0}})$$
for all $t\in\rz$, then
$$\e^{tr}\EH\in\big(\rqom{}\times\dqpeom{}\big)\cap\qh{2}{q,q+1}{}{}(\Bc_{r_{0}+1})$$
holds for all $t\in\rz$.
The assertion holds in particular if $\FG$ is compactly supported.
\end{theo}

\begin{proof}
The idea of the proof is to estimate the exponential series. 
For this, we need some technical preliminaries.
For all $s\in\rz$ we have
$$(M+\ie\omega\Lambda)\EH=-\ie\Lambda\FG\in\Lzqqpesom\cap\qh{2}{q,q+1}{s}{}(\Bc_{r_{0}}).$$
Hence, by Theorem \ref{polyabklmax} and Remark \ref{bempolyabklmax} 
$\EH$ belongs to $\rqsom\times\dqpesom$ for all $s\in\rz$. 
Therefore, inner regularity, e.g. a combination of a cutting technique together with \cite[Lemma 4.2(i)]{kuhnpaulyreg},
yields $\EH\in\qh{2}{q,q+1}{s}{}(\Bc_{r_{0}+\theta})$ for all $s\in\rz$ and all positive $\theta$. Thus, 
$$(M+\ie\omega)\EH=-\ie\Lambda\FG-\ie\omega\Lambdad\EH=:\FGs\in\Lzqqpesom\cap\qh{2}{q,q+1}{s}{}(\Bc_{r_{0}+\theta})$$
for all $\theta>0$ and all $s\in\rz$. Consequently, applying $\pdiv$, $\rot$ and $(M-\ie\omega)$ 
to the latter equation we can compute in $\Bc_{r_{0}}$
\beq\ie\omega\pdiv E=\pdiv\Fs,\quad\ie\omega\rot H=\rot\Gs\mylabel{deltadEH}\eeq
and
$$(M^2+\omega^2)\EH=(M-\ie\omega)\FGs.$$
By defining $\square:=\Delta-M^2$ we note
$$M^2=\zmat{\pdiv\rot}{0}{0}{\rot\pdiv},\quad\square=\zmat{\rot\pdiv}{0}{0}{\pdiv\rot}
,\quad M^2+\square=\Delta.$$
Therefore, we achieve
$$(\Delta+\omega^2)\EH=(M-\ie\omega-\frac{\ie}{\omega}\square)\FGs=:\FGd,$$
which is the equation we intend to work with. 
Now, we multiply this equation and all forms by $r^m$ with some $m\in\rz$
and indicate the resulting forms by an index $m$.
We note that all occurring forms are well defined elements of
$\qh{2}{q,q+1}{}{}(\Bc_{r_{0}+\theta})$ for all $\theta>0$ and all $m\in\rz$.
Using Lemma \ref{gewichttausch} we obtain in $\Bc_{r_{0}}$
\begin{align}
\begin{split}
&\qquad\big(\Delta+\omega^2+\frac{m(m+2-N)}{r^2}\big)\eh-2\frac{m}{r}\p_r\eh\\
&=\fgd=(M-\ie\omega-\frac{\ie}{\omega}\Box)\fgs+(-C_{M,r^m}+\frac{\ie}{\omega}C_{\Box,r^m})\FGs\\
&=-\ie\big(\ie M+\omega+\frac{1}{\omega}\Box-\frac{1}{\omega}\frac{m}{r}
(\ie\omega S+\sw MS+\sw SM-\frac{m+1}{r}\sw S^2)\big)\fgs.
\end{split}\mylabel{gewgleich}
\end{align}
With $\eta $ from Lemma \ref{decomplemma} we define the cut-off function 
$\varphi_{\theta}:=\eta(r-\theta+1)$ for all $\theta>r_{0}+1$. Then, 
$$\supp\varphi_{\theta}\subset\Bc_{\theta},\quad 
\supp\nabla\varphi_{\theta},\supp(1-\varphi_{\theta})\cap\Bc_{\theta}
\subset\Bc_{\theta}\cap B_{\theta+1}.$$
Without loss of generality let any form be real-valued.
We multiply \eqref{gewgleich} by $\varphi_\theta r^p\eh$ with $p\in\rz$ resp. 
$\varphi_\theta r\p_r\eh$ and integrate over $\rN$. We achieve
\begin{align}
\begin{split}
&\qquad\skp{\fgd}{\varphi_\theta r^p\eh}_{\Lzqqpe{}}\\
&=\skpb{(\Delta+\omega^2+\frac{m(m+2-N)}{r^2})\eh-2\frac{m}{r}\p_r\eh}{\varphi_\theta r^p\eh}_{\Lzqqpe{}}\\
&=\int_{\rN}\Delta\eh\varphi_\theta r^p\eh
+\omega^2\varphi_\theta r^p|\eh|_{q,q+1}^2\\
&\qquad+m(m+2-N)\varphi_\theta r^{p-2}|\eh|_{q,q+1}^2
-m\varphi_\theta r^{p-1}\p_{r}|\eh|_{q,q+1}^2\dla
\end{split}\mylabel{fgsmulte}
\intertext{resp.}
\begin{split}
&\qquad\skp{\fgd}{\varphi_\theta r\p_r\eh}_{\Lzqqpe{}}\\
&=\skpb{(\Delta+\omega^2+\frac{m(m+2-N)}{r^2})\eh-2\frac{m}{r}\p_r\eh}{\varphi_\theta r\p_r\eh}_{\Lzqqpe{}}\\
&=\int_{\rN}\Delta\eh\varphi_\theta r\p_{r}\eh
+\frac{\omega^2}{2}\varphi_\theta r\p_{r}|\eh|_{q,q+1}^2\\
&\qquad+\frac{m}{2}\frac{m+2-N}{r}\varphi_\theta\p_{r}|\eh|_{q,q+1}^2
-2m\varphi_\theta|\p_{r}\eh|_{q,q+1}^2\dla.
\end{split}\mylabel{fgsmultz}
\end{align}
By partial integration we get from \eqref{fgsmulte} resp. \eqref{fgsmultz} 
\begin{align}
\begin{split}
&\qquad\Big|\int_{\rN}\varphi_\theta r^p\big((\omega^2+\frac{\gamma_{m,p,N}}{r^2})|\eh|_{q,q+1}^2
-\sum_{|\alpha|=1}|\pa\eh|_{q,q+1}^2\big)\dla\Big|\\
&\leq\big|\skp{\fgd}{\varphi_\theta r^p\eh}_{\Lzqqpe{}}\big|+c m(\theta+1)^{2m+p}
\end{split}\mylabel{schmutzone}
\intertext{resp.}
\begin{split}
&\qquad\Big|\int_{\rN}\varphi_\theta\big((N\omega^2+\frac{\tilde{\gamma}_{m,N}}{r^2})|\eh|_{q,q+1}^2
-(N-2)\sum_{|\alpha|=1}|\pa\eh|_{q,q+1}^2\\
&\qquad\qquad+4m|\p_r\eh|_{q,q+1}^2\big)\dla\Big|\\
&\leq\big|\skp{\fgd}{\varphi_\theta r\p_r\eh}_{\Lzqqpe{}}\big|+c m^2(\theta+1)^{2m+1},
\end{split}\mylabel{schmutztwo}
\end{align}
where 
$$\gamma_{m,p,N}:=m(m+p)+p(p+N-2)/2,\quad \tilde{\gamma}_{m,N}:=m(m+2-N)(N-2)$$
and $c$ is a generic constant independent of $m$ and $\theta$.
Now, we multiply \eqref{schmutzone} for $p=0$ by $N-2$ and add
\eqref{schmutzone} and \eqref{schmutztwo} in a suitable way. We obtain
\begin{align}
&\qquad\int_{\rN}\varphi_\theta\big((2\omega^2-\frac{m(N-2)^2}{r^2})|\eh|_{q,q+1}^2
+4m|\p_r\eh|_{q,q+1}^2\big)\dla\non\\
&\leq c\Big(\big|\skp{\fgd}{\varphi_\theta \eh}_{\Lzqqpe{}}\big|
+\big|\skp{\fgd}{\varphi_\theta r\p_r\eh}_{\Lzqqpe{}}\big|
+m^2(\theta+1)^{2m+1}\Big),\non
\end{align}
which yields
\begin{align}
\begin{split}
\int_{\rN}\varphi_\theta|\eh|_{q,q+1}^2\dla
&\leq(N-2)^2\int_{\rN}\varphi_\theta\frac{1}{r}\frac{m}{r}|\eh|_{q,q+1}^2\dla+cm^2(\theta+1)^{2m+1}\\
&\qquad+c\Big(\big|\skp{\fgd}{\varphi_\theta \eh}_{\Lzqqpe{}}\big|
+\big|\skp{\fgd}{\varphi_\theta r\p_r\eh}_{\Lzqqpe{}}\big|\Big).
\end{split}\mylabel{ehabsche}
\end{align}
To take care of the right hand side we prove two lemmas.

\begin{lem}\mylabel{zweiteablabsch}
For all $p\in\rz$ there exists a constant $c>0$, such that for all $\psi\in\h{2}{p/2}{}(\Bc_{r_{0}+1})$
and all $\theta>r_{0}+1$
$$\sum_{|\alpha|=2}\int_{\rN}\varphi_\theta r^p|\pa\psi|_{q,q+1}^2\dla
\leq c\int_{\Bc(\theta)}r^p\big(\sum_{|\alpha|=1}|\pa\psi|_{q,q+1}^2+|\Delta\psi|_{q,q+1}^2\big)\dla.$$
\end{lem}

\begin{proof}
By a cutting-technique we may assume $\psi\in\h{2}{p/2}{\circ}(\Bc_{r_{0}+1})$.
Hence, it suffices to prove the assertion only for $\psi\in\cqun(\Bc_{r_{0}+1})$ by continuity.
But for those $\psi$ the result is simply shown by several partial integrations.
\end{proof}
 
\begin{lem}\mylabel{schmutzlem}
Let $p,\tilde{p},\hat{p}\in\rz$ with $\tilde{p}>p$ and $\hat{p}>2$.
There are a constant $c>0$ and a non-negative function $\varkappa$ tending to zero at infinity,
such that for all $m\in\rz$ and $\sigma>0$, $\theta>r_{0}+1$ the following estimates hold:
\begin{align*}
\text{\rm\bf(i)}&&
\big|\skp{\fgd}{\varphi_\theta r^p\eh}_{\Lzqqpe{}}\big|
&\leq\varkappa(\theta)\Big(m(\theta+1)^{2m+p}
+\int_{\rN}\varphi_\theta r^{\tilde{p}}|\fg|_{q,q+1}^2\dla\\
&&&\qquad+\int_{\rN}\varphi_\theta r^p\big(\frac{m^4}{r^4}|\eh|_{q,q+1}^2
+\sum_{|\alpha|\leq2}|\pa\eh|_{q,q+1}^2\big)\dla\Big)\\
\text{\rm\bf(ii)}&&
\big|\skp{\fgd}{\varphi_\theta r^p\eh}_{\Lzqqpe{}}\big|
&\leq\varkappa(\theta)\Big(m(\theta+1)^{2m+p}
+\sum_{|\alpha|\leq1}\int_{\rN}\varphi_\theta r^{\tilde{p}}|\pa\fg|_{q,q+1}^2\dla\\
&&&\qquad+\int_{\rN}\varphi_\theta r^p\big(\frac{m^4}{r^4}|\eh|_{q,q+1}^2
+\sum_{|\alpha|\leq1}|\pa\eh|_{q,q+1}^2\big)\dla\Big)\\
&&&\qquad+c\int_{\rN}\varphi_\theta r^p\big(\frac{1}{\sigma}|\eh|_{q,q+1}^2
+\sigma\sum_{|\alpha|=1}|\pa\eh|_{q,q+1}^2\big)\dla\\
\text{\rm\bf(iii)}&&
\big|\skp{\fgd}{\varphi_\theta r\p_r\eh}_{\Lzqqpe{}}\big|
&\leq\varkappa(\theta)\Big(m^2(\theta+1)^{2m+1}
+\sum_{|\alpha|\leq2}\int_{\rN}\varphi_\theta r^{\hat{p}}|\pa\fg|_{q,q+1}^2\dla\\
&&&\qquad+\int_{\rN}\varphi_\theta\big(\sum_{|\alpha|\leq1}\frac{m^4}{r^4}|\pa\eh|_{q,q+1}^2
+\sum_{|\alpha|\leq2}|\pa\eh|_{q,q+1}^2\big)\dla\Big)
\end{align*}
\end{lem}

\begin{proof}
By several partial integrations we remove all derivatives from $\Lambdad\eh$ resp. $\fg$ 
and use the decay of $\Lambdad$ resp. the integrability of $\fg$. This yields {\bf(i)}.
In {\bf(ii)} there are no longer second derivatives allowed on $\eh$.
Hence, we have to insert a $\sigma$ into the estimate 
since we do not demand any decay properties of the derivatives of $\Lambdad$.
To prove {\bf(iii)} we have to handle one challenging term, i.e.
\begin{align*}
&\qquad\skp{\Box\Lambdad\eh}{\varphi_\rho r\p_r\eh}_{\Lzqqpe{}}\\
&=-\skpb{\pdiv\epsd E_{m}}{\pdiv(\varphi_\rho r\p_r E_{m})}_{\qLz{q-1}{}}
-\skpb{\rot\mud H_{m}}{\rot(\varphi_\rho r\p_r H_{m})}_{\qLz{q+2}{}}.
\end{align*}
The part of the first summand causing the biggest difficulties is
$$\skp{\pdiv\epsd E_{m}}{\varphi_\rho r\p_r\pdiv E_{m}}_{\qLz{q-1}{}}$$
because we can not integrate by parts anymore since $E$ is only twice weakly differentiable.
But with \eqref{deltadEH} and Lemma \ref{gewichttausch} (ii) we can substitute
$$\pdiv E_{m}=\frac{m}{r}TE_{m}
-\frac{1}{\omega}(\pdiv-\frac{m}{r}T)\eps F_{m}
-(\pdiv-\frac{m}{r}T)\epsd E_{m}.$$
Thus, the most challenging term reads now
$$\skp{\pdiv\epsd E_{m}}{\varphi_\rho r\p_r\pdiv\epsd E_{m}}_{\qLz{q-1}{}}
=\peh\int_{\rN}\varphi_\rho r\p_r|\pdiv\epsd E_{m}|_{q-1}^2\dla$$
and can be handled easily by two partial integrations.
\end{proof}

We proceed with the proof of the theorem by further estimating \eqref{ehabsche}
using Lemma \ref{schmutzlem} (i) and (iii).
\begin{align}
\begin{split}
&\qquad\int_{\rN}\varphi_\theta |\eh|_{q,q+1}^2\dla\\
&\leq\varkappa(\theta)\big(\sum_{|\alpha|\leq1}\int_{\rN}\varphi_\theta(1+\frac{m^4}{r^4})|\pa\eh|_{q,q+1}^2\dla
+\sum_{|\alpha|=2}\int_{\rN}\varphi_\theta |\pa\eh|_{q,q+1}^2\dla\\
&\qquad+\sum_{|\alpha|\leq2}\int_{\rN}\varphi_\theta r^{\hat{p}}|\pa\fg|_{q,q+1}^2\dla\big)
+cm^2(\theta+1)^{2m+1}
\end{split}\mylabel{ehabschone}
\intertext{Utilizing Lemma \ref{zweiteablabsch} we estimate the second derivatives of $\eh$ 
by the first ones and $\Delta\eh$. Then we substitute $\Delta\eh$ with \eqref{gewgleich}
and get for sufficient large $\theta$}
&\qquad\sum_{|\alpha|=2}\int_{\rN}\varphi_\theta|\pa\eh|_{q,q+1}^2\dla\non\\
&\leq c\big(\sum_{|\alpha|\leq1}\int_{\rN}\varphi_\theta(1+\frac{m^4}{r^4})|\pa\eh|_{q,q+1}^2\dla\non\\
&\qquad+\sum_{|\alpha|\leq2}\int_{\rN}\varphi_\theta(1+\frac{m^4}{r^4})|\pa\fg|_{q,q+1}^2\dla\big)
+cm^4(\theta+1)^{2m}.\non
\intertext{Now we insert this estimate into \eqref{ehabschone} and obtain}
\begin{split}
&\qquad\int_{\rN}\varphi_\theta |\eh|_{q,q+1}^2\dla\\
&\leq\varkappa(\theta)\big(\sum_{|\alpha|\leq1}\int_{\rN}\varphi_\theta(1+\frac{m^4}{r^4})|\pa\eh|_{q,q+1}^2\dla\\
&\qquad+\sum_{|\alpha|\leq2}\int_{\rN}\varphi_\theta r^{\hat{p}}(1+\frac{m^4}{r^4})|\pa\fg|_{q,q+1}^2\dla\big)
+cm^4(\theta+1)^{2m+1}.
\end{split}\mylabel{ehabschz}
\intertext{For all $p\in\rz$ \eqref{schmutzone} yields the estimate}
&\qquad\sum_{|\alpha|=1}\int_{\rN}\varphi_\theta r^p|\pa\eh|_{q,q+1}^2\dla\non\\
&\leq c\int_{\rN}\varphi_\theta r^p(1+\frac{m^2}{r^2})|\eh|_{q,q+1}^2\dla
+\big|\skp{\fgd}{\varphi_\theta r^p\eh}_{\Lzqqpe{}}\big|+c m(\theta+1)^{2m+p},\non
\intertext{such that by Lemma \ref{schmutzlem} (ii) for sufficient small $\sigma$ and large $\theta$}
&\qquad\sum_{|\alpha|=1}\int_{\rN}\varphi_\theta r^p|\pa\eh|_{q,q+1}^2\dla\non\\
&\leq c \int_{\rN}\varphi_\theta r^p(1+\frac{m^4}{r^4})|\eh|_{q,q+1}^2\dla\non\\
&\qquad+c\sum_{|\alpha|\leq1}\int_{\rN}\varphi_\theta r^{\tilde{p}}|\pa\fg|_{q,q+1}^2\dla
+cm(\theta+1)^{2m+p}\non
\intertext{follows. Now, we plug the latter estimate for $p=0$ and $p=-4$ into \eqref{ehabschz} obtaining}
&\qquad\int_{\rN}\varphi_\theta |\eh|_{q,q+1}^2\dla\non\\
&\leq\varkappa(\theta)\big(\int_{\rN}\varphi_\theta(1+\frac{m^8}{r^8})|\eh|_{q,q+1}^2\dla\non\\
&\qquad+\sum_{|\alpha|\leq2}\int_{\rN}\varphi_\theta m^4r^{\hat{p}}|\pa\fg|_{q,q+1}^2\dla\big)
+cm^5(\theta+1)^{2m+1},\non
\intertext{where we assume without loss of generality $\tilde{p}\leq{\hat{p}}$.
Therefore, for sufficient large $\theta$}
&\qquad\int_{\rN}\varphi_\theta |\eh|_{q,q+1}^2\dla\non\\
&\leq\varkappa(\theta)\big(\int_{\rN}\varphi_\theta \frac{m^8}{r^8}|\eh|_{q,q+1}^2\dla\non\\
&\qquad+\sum_{|\alpha|\leq2}\int_{\rN}\varphi_\theta m^4r^{\hat{p}}|\pa\fg|_{q,q+1}^2\dla\big)
+cm^5(\theta+1)^{2m+1},\non
\intertext{i.e.}
&\qquad\int_{\Bc_{\theta}}|\eh|_{q,q+1}^2\dla\non\\
&\leq\varkappa(\theta)\big(\int_{\Bc_{\theta}}\frac{m^8}{r^8}|\eh|_{q,q+1}^2\dla
+\sum_{|\alpha|\leq2}\int_{\Bc_{\theta}}m^4r^{\hat{p}}|\pa\fg|_{q,q+1}^2\dla\big)
+cm^5(\theta+1)^{2m+1}.\non
\end{align}
Setting $k:=2m$ we finally get for all $k\in\rz$ and sufficient large $\theta$
\begin{align}
\begin{split}
&\qquad\int_{\Bc_{\theta}}r^k|\EH|_{q,q+1}^2\dla\\
&\leq\varkappa(\theta)k^8\big(\int_{\Bc_{\theta}}r^{k-8}|\EH|_{q,q+1}^2\dla
+\sum_{|\alpha|\leq2}\int_{\Bc_{\theta}}r^{k+{\hat{p}}}|\pa\FG|_{q,q+1}^2\dla\big)
+ck^5(\theta+1)^{k+1}.
\end{split}\mylabel{expabklabsch}
\end{align}

We are ready to prove the exponential decay.
Let $t\in\rzp$ and $\theta$ large enough, such that \eqref{expabklabsch} is holds.
Then, we have for all natural numbers $8\leq K_1\leq K_2$
\begin{align*}
&\qquad\sum_{k=K_1}^{K_2}\frac{t^k}{k!}\int_{\Bc_{\theta}}r^k|\EH|_{q,q+1}^2\dla\\
&\leq\varkappa(\theta)\sum_{k=K_1}^{K_2}\frac{t^k}{k!}k^8\int_{\Bc_{\theta}}r^{k-8}|\EH|_{q,q+1}^2\dla\\
&\qquad+\varkappa(\theta)\sum_{k=K_1}^{K_2}\frac{t^k}{k!}\sum_{|\alpha|\leq2}\int_{\Bc_{\theta}}\ub{k^8r^{k+{\hat{p}}}}_{\leq8^kr^kr^{\hat{p}}}|\pa\FG|_{q,q+1}^2\dla
+c\sum_{k=K_1}^{K_2}\frac{t^k}{k!}\ub{k^5(\theta+1)^{k+1}}_{\leq(5(\theta+1))^k(\theta+1)}\\
&\leq\varkappa(\theta)\sum_{k=K_1}^{K_2}\frac{t^k}{(k-8)!}\int_{\Bc_{\theta}}r^{k-8}|\EH|_{q,q+1}^2\dla\\
&\qquad+\varkappa(\theta)\sum_{|\alpha|\leq2}\int_{\Bc_{\theta}}\e^{8tr}r^{\hat{p}}|\pa\FG|_{q,q+1}^2\dla
+c\e^{5t(\theta+1)}(\theta+1)\\
&\leq\varkappa(\theta)\sum_{k=K_1-8}^{K_2}\frac{t^{k+8}}{k!}\int_{\Bc_{\theta}}r^k|\EH|_{q,q+1}^2\dla
+\varkappa(\theta)+c\e^{6t(\theta+1)}.
\end{align*}
Now, let $\theta$ be so large, such that $\varkappa(\theta)t^8\leq1/2$. Then,
$$\sum_{k=K_1}^{K_2}\frac{t^k}{k!}\int_{\Bc_{\theta}}r^k|\EH|_{q,q+1}^2\dla
\leq\varkappa(\theta)\sum_{k=0}^{K_1-1}\frac{t^{k+8}}{k!}\int_{\Bc_{\theta}}r^k|\EH|_{q,q+1}^2\dla+c.$$
Since the right hand side is independent of $K_2$ 
we obtain by the monotone convergence theorem for $K_{2}\to\infty$
$$\int_{\Bc_{\theta}}\e^{tr}|\EH|_{q,q+1}^2\dla<\infty,$$
i.e. $\e^{tr}\EH\in\Lzqqpeom{}$ for all $t\in\rz$. 
The differential equation yields 
$$\e^{tr}M\EH\in\Lzqqpeom{},\quad\e^{tr}(\pdiv\eps E,\rot\mu H)\in\qLz{q-1,q+1}{}(\Bc_{r_{0}+\theta})$$
for all $t\in\rz$ and all $\theta>0$. Consequently, by
\begin{align*}
M\e^{tr}\EH&=\e^{tr}M\EH+t\e^{tr}S\EH,\\
\pdiv\eps\e^{tr}E&=\e^{tr}\pdiv\eps E+t\e^{tr}T\eps E,\quad
\rot\mu\e^{tr}H=\e^{tr}\rot\mu H+t\e^{tr}R\mu H
\end{align*}
and inner regularity (see the beginning of the proof) we achieve
$\e^{tr}\EH\in\qh{1}{q,q+1}{}{}(\Bc_{r_{0}+\theta})$
for all $t\in\rz$ and all $\theta>0$. Repeating this argument yields finally
the same assertion for $\h{2}{}{}$.
\end{proof}

\section{Solution theory for time-harmonic Maxwell equations}

As a canonical application we intend to present a solution theory 
for the radiation problem \eqref{shortproblem}.
We will follow in close lines the first part of \cite{xmas}.
It can be seen easily that
$$\rot:\ronqom{}\subset\lzqom\to\lzqpeom,\quad\pdiv:\dqpeom{}\subset\lzqpeom\to\lzqom$$
are skew-adjoint to each other. Consequently, for $0$-admissible coefficients $\Lambda$ 
$$\calM:D(\calM):=\ronqom{}\times\dqpeom{}\subset{}_\Lambda\Lzqqpeom{}\to{}_\Lambda\Lzqqpeom{},\quad
\calM\EH:=M_{\Lambda}\EH$$
is self-adjoint.
Here, ${}_\Lambda\Lzqqpeom{}:=\Lzqqpeom{}$ equipped with the scalar product 
$$(\EH,v)\mapsto\skp{\Lambda u}{v}_{\Lzqqpeom{}}.$$
Furthermore, we will denote the kernel and the image of $\calM$ 
by $N(\calM)$ and $I(\calM)$, respectively.

\begin{defini}\mylabel{defloesc}
Let $\omega\in\cz\ohne\rz$ and $\FG\in\Lzqqpelocom$. 
Then $\EH$ solves $\Max(\Lambda,\omega,\FG)$, if and only if 
$\EH\in D(\calM)$ and $(M_{\Lambda}-\omega)\EH=\FG$.
\end{defini}

The self-adjointness of $\calM$ yields the unique solvability of $\Max(\Lambda,\omega,\FG)$
for non-real frequencies $\omega\in\cz\ohne\rz$ and right hand sides $\FG\in\Lzqqpeom{}$.
We denote the continuous solution operator by $\loesom$.
Since the spectrum of $\calM$ is contained in the real axis
we expect from well known facts about Helmholtz' equation 
that we have to work in weighted $\lz$-spaces 
and utilize radiating solutions to get a proper solution theory for real frequencies.

\begin{defini}\mylabel{defloesr}
Let $\omega\in\rz\ohne\{0\}$ and $\FG\in\Lzqqpelocom$. 
Then $\EH$ solves $\Max(\Lambda,\omega,\FG)$, if and only if
$$\EH\in\ronqkmehom\times\dqpekmehom,\quad
(S+1)\EH\in\Lzqqpegmehom,\quad
(M_{\Lambda}-\omega)\EH=\FG.$$
The second constraint will be called
'Maxwell incoming radiation condition' or simply 'radiation condition'.
\end{defini}

We will establish a solution theory using Eidus' limiting absorption principle.
The key tool for the application of this principle is an a priori estimate,
which ensures the uniform continuity of $\loesom$ operating in proper Hilbert spaces 
even up to the real axis.

\subsection{An a priori estimate}

\begin{lem}\mylabel{apriorimaxwell}
Let $I\subset\rzon$ be a compact interval and $s,-t>1/2$. 
Furthermore, let $\Lambda$ be $\tau$-admissible with $\tau>1$.
Then, there exist constants $c,\theta>0$ and a $\hat{t}>-1/2$, 
such that for all $\omega\in\cz_+$,
which means that $\omega$ has got non negative imaginary part,
with $\omega^2=\lambda^2+\ie\sigma\lambda$,
$\lambda\in I$, $\sigma\in(0,1]$ and $\FG\in\Lzqqpesom$
as well as $\EH:=\loesom\FG$ the estimate
$$\norm{\EH}_{\rqtom\times\dqpetom}
+\normb{(S+1)\EH}_{\Lzqqpeom{\hat{t}}}
\leq c\big(\norm{\FG}_{\Lzqqpesom}
+\norm{\EH}_{\Lzqqpe{}(\Omega\cap B_{\theta})}\big)$$
holds true.
\end{lem}

\begin{proof}
Without loss of generality let $s\in(1/2,1)$. We note $(1/2,1)\cap\pI=\emptyset$.
Decomposing $\FG$ and $\EH$ using Lemma \ref{decomplemma} 
with $s=s$, $t=0$ we have $\FSGS=0$ since $s<N/2$
and obtain $\EHD\in\qh{2}{q,q+1}{}{}$ satisfying
$$(\Delta+\omega^2)\EHD=-(1+\omega^2)\EHF+(1-\ie\omega)\FGs=:\FGD\in\Lzqqpes.$$
The self-adjointness of $\Delta:\qh{2}{q,q+1}{}{}\subset\Lzqqpe{}\to\Lzqqpe{}$ yields
$(\Delta+\omega^2)^\me\FGD=\EHD$.
Applying \cite[Lemma 7]{linelae}, which is a well known a priori estimate 
for the scalar Helmholtz equation in $\rN$;
see also Ikebe and Saito \cite{ikebesaito} or Vogelsang \cite[section 2]{volker},
componentwise to $\EHD$ and using Lemma \ref{decomplemma} with
$M(\e^{-\ie\lambda r}\EHD)=\e^{-\ie\lambda r}(M-\ie\lambda S)\EHD$
we get the estimate
\begin{align}
\begin{split}
&\qquad\norm{\EHD}_{\Lzqqpet}+\normb{(M-\ie\lambda S)\EHD}_{\Lzqqpe{s-1}}\\
&\leq c\Big(\norm{\EHD}_{\Lzqqpet}
+\sum_{|\alpha|=1}\normb{\pa(\e^{-\ie\lambda r}\EHD)}_{\Lzqqpe{s-1}}\Big)\\
&\leq c\norm{\FGD}_{\Lzqqpes}
\leq c\big(\norm{\FG}_{\Lzqqpesom}+\norm{\EH}_{\Lzqqpeom{s-\tau}}\big),
\end{split}\mylabel{aprioriabschlambda}
\end{align}
which holds uniformly in $\EHD$, $\FGD$ and $\lambda$.
But actually we are interested in estimating the term 
$\normb{(M-\ie\omega S)\EHD}_{\Lzqqpe{s-1}}$.
This needs an additional argument. The standard resolvent estimate yields
\beq\sigma|\lambda|\norm{\EHD}_{\Lzqqpe{}}
\leq\norm{\FGD}_{\Lzqqpe{}}\mylabel{sigmatrick}\eeq
and we have $|\omega+\lambda|\geq|\lambda|$
since $|\ret\omega|\geq|\lambda|\sqrt2/2$ and $\omega\in\czp$.
Thus, by \eqref{sigmatrick} and
$$\omega-\lambda=\frac{\omega^2-\lambda^2}{\omega+\lambda}=\frac{\ie\sigma\lambda}{\omega+\lambda}$$
we achieve uniformly in $\omega$
\begin{align*}
\normb{(M-\ie\omega S)\EHD}_{\Lzqqpe{s-1}}
&\leq\normb{(M-\ie\lambda S)\EHD}_{\Lzqqpe{s-1}}+c|\omega-\lambda|\norm{\EHD}_{\Lzqqpe{s-1}}\\
&\leq\normb{(M-\ie\lambda S)\EHD}_{\Lzqqpe{s-1}}+\frac{c}{|\lambda|}\norm{\FGD}_{\Lzqqpe{}}.
\end{align*}
A combination of the latter estimate with \eqref{aprioriabschlambda} and Lemma \ref{decomplemma} yield
\begin{align*}
&\qquad\norm{\EH}_{\rqtom\times\dqpetom}+\normb{(M-\ie\omega S)\EH}_{\Lzqqpeom{s-1}}\\
&\leq c\Big(\norm{\EHD}_{\Lzqqpet}+\normb{(M-\ie\omega S)\EHD}_{\Lzqqpe{s-1}}
+\norm{\FG}_{\Lzqqpesom}+\norm{\EH}_{\Lzqqpeom{s-\tau}}\Big)\\
&\leq c\big(\norm{\FG}_{\Lzqqpesom}+\norm{\EH}_{\Lzqqpeom{s-\tau}}\big)
\end{align*}
uniformly in $\EH$, $\FG$ and $\omega$. 
Using $(M-\ie\omega S)\EH=-\ie\omega(S+1+\Lambdad)\EH+\FG$ we finally arrive at
$$\norm{\EH}_{\rqtom\times\dqpetom}+\normb{(S+1)\EH}_{\Lzqqpeom{s-1}}
\leq c\big(\norm{\FG}_{\Lzqqpesom}+\norm{\EH}_{\Lzqqpeom{s-\tau}}\big).$$
Because of the monotone dependence of the weighted $\lz$-norms on the weights we may assume 
$t$ near to $-1/2$ and $s$ near to $1/2$, such that $1<s-t<\tau$ holds. 
Then Lemma \ref{normtrick} completes the proof.
\end{proof}

\subsection{Fredholm theory using the limiting absorption principle}

We prove three more technical lemmas.

\begin{lem}\mylabel{techstrahl}
Let $\alpha,\beta\in\rz$ with $0\leq\alpha<\beta$ and $\rN\ohne\Omega\subset B_{\alpha}$.
Moreover, for some $t\in\rz$ let $E\in\ronqtom$, $H\in\dqpetom$
and $\varphi\in\pc{0}\big([\alpha,\beta],\cz\big)$.
Then with
$$\Abb{\psi}{[0,\beta]}{\cz}{\sigma}{\bds\int_{\max\{\alpha,\sigma\}}^{\beta}\varphi(s)\,ds\eds}$$
and $\Phi:=\varphi\circ r$, $\Psi:=\psi\circ r$
$$\skp{\Phi RE}{H}_{\lzqpe(\Bc_{\alpha}\cap B_{\beta})}
=\skp{\Psi\rot E}{H}_{\lzqpe(\Omega\cap B_{\beta})}
+\skp{\Psi E}{\pdiv H}_{\lzq(\Omega\cap B_{\beta})}.$$
\end{lem}

\begin{proof}
Assume $E\in\cqunom$, $H\in\cqpeuom$. 
With $\gamma:=\bds\int_{\alpha}^{\beta}\varphi(s)\,ds\eds$ we have
$\restr{\psi}{[0,\alpha]}=\gamma$, $\psi(\beta)=0$ 
and $\psi\in\pc{1}\big((\alpha,\beta),\cz\big)$ with $\psi'=-\varphi$.
By Stokes' theorem we compute
\begin{align*}
&\qquad\skp{\Psi\rot E}{H}_{\lzqpe(\Omega\cap B_{\beta})}
+\skp{\Psi E}{\pdiv H}_{\lzq(\Omega\cap B_{\beta})}\\
&=\gamma\skp{\rot E}{H}_{\lzqpe(\Omega\cap B_{\alpha})}
+\gamma\skp{E}{\pdiv H}_{\lzq(\Omega\cap B_{\alpha})}\\
&\qquad+\skp{\Psi\rot E}{H}_{\lzqpe(\Bc_{\alpha}\cap B_{\beta})}
+\skp{\Psi E}{\pdiv H}_{\lzq(\Bc_{\alpha}\cap B_{\beta})}\\
&=\gamma\int_{S_{\alpha}}\iota_{\alpha}^*(E\wedge*\bar{H})
+\skp{\Phi RE}{H}_{\lzqpe(\Bc_{\alpha}\cap B_{\beta})}\\
&\qquad-\gamma\int_{S_{\alpha}}\iota_{\alpha}^*(E\wedge*\bar{H})
+\psi(\beta)\int_{S_{\beta}}\iota_{\beta}^*(E\wedge*\ol{H}),
\end{align*}
where $\iota_\theta:S_\theta\to\rN$ denotes the natural embedding.
With the help of mollifiers we get the desired formula for all $H\in\dqpetom$.
Since $\cqunom$ is dense in $\ronqtom$ the assertion holds as stated.
\end{proof}

By the same approximation technique we obtain the rule of partial integration for weighted forms.

\begin{lem}\mylabel{schaufeln}
Let $t,s\in\rz$ and $E\in\ronqtom$, $H\in\dqpesom$
as well as $\varphi_\theta:=1-\eta(\,\cdot\,/\theta)$,
$\Phi_\theta:=\varphi_\theta\circ r$ with some $\theta>r_{0}$. Then
$$\skp{\Phi_\theta \rot E}{H}_{\lzqpeom}+\skp{\Phi_\theta E}{\pdiv H}_{\lzqom}=-\skpb{\varphi_\theta'(r)RE}{H}_{\lzqpeom}$$
holds. Additionally, if $t+s\geq0$ then
$$\skp{\rot E}{H}_{\lzqpeom}+\skp{E}{\pdiv H}_{\lzqom}=0.$$
\end{lem}

\begin{rem}\mylabel{schaufelnbem}
If $\EH\in\ronqtom\times\dqpetom$ and
$v\in\ronqsom\times\dqpesom$ with $t+s\geq0$ then
$$\skp{M\EH}{v}_{\Lzqqpeom{}}+\skp{\EH}{Mv}_{\Lzqqpeom{}}=0.$$
\end{rem}

Since $\calM$ is self-adjoint we have the well known Hodge-Helmholtz decomposition.

\begin{lem}\mylabel{hodgehelmholtzdeco}
\begin{align*}
\Lzqqpeom{}={}_{\Lambda}\Lzqqpeom{}
&=N(\calM)\oplus_{\Lambda}\ol{I(\calM)}\\
&=\big(\ronqnom{}\times\dqpenom{}\big)\oplus_{\Lambda}
\Lambda^\me\big(\ol{\pdiv\dqpeom{}}\times\ol{\rot\ronqom{}}\big)
\end{align*}
Here $\oplus_{\Lambda}$ denotes the orthogonal sum in ${}_{\Lambda}\Lzqqpeom{}$
and the closures are taken in the respective $\lz$-spaces.
\end{lem}

Another essential ingredient of the solution theory 
generating convergence in the limiting absorption argument
is the so called Maxwell local compactness property {\sf MLCP}, 
i.e. the embeddings
$$\ronqsom\cap\dqsom\hookrightarrow\Lzqtom$$
have to be compact for all $t<s$ and all $q$.
For a detailed analysis of this property of $\dom$
we refer to \cite{weckmax}, \cite{xmas}
and \cite{paulytimeharm,paulystatic,paulydeco,paulyasym} 
as well as the papers cited there.

\begin{defini}\mylabel{punktspektrum}
For $\omega\in\cz\ohne\{0\}$ we define
\begin{align*}
P&:=\setb{\omega\in\cz\ohne\{0\}}{\Max(\Lambda,\omega,0)\text{\rm\, has a nontrivial solution.}},\\
N_{\omega}&:=\setb{\EH}{\EH\text{\rm\, is a solution of }\Max(\Lambda,\omega,0).}.
\end{align*}
\end{defini}

We remark $P\subset\rzon$ or $N_{\omega}=N(\calM-\omega)=\{0\}$ if $\omega\in\cz\ohne\rz$.
We are ready to prove the main result of this section.

\begin{theo}\mylabel{fredholm}
Let $\omega\in\rzon$ and $\Lambda$ be $\tau$-admissible with $\tau>1$.
\begin{itemize}
\item[\bf(i)] Eigen-solutions decay polynomially, i.e. for all $t\in\rz$
\begin{align*}
N_{\omega}&=N(\calM-\omega)
\subset\big(\ronqtom\cap\eps^\me\dqtnom\big)\times\big(\dqpetom\cap\mu^\me\ronqpetnom\big).
\end{align*}
\end{itemize}
Additionally, let $\Omega$ possess the {\sf MLCP}. Then
\begin{itemize}
\item[\bf(ii)] $N_{\omega}$ is finite dimensional;
\item[\bf(iii)] $P$ has no accumulation point in $\rzon$;
\item[\bf(iv)] for every $\FG\in\Lzqqpegehom$ there exists a solution $\EH$ of the problem
$\Max(\Lambda,\omega,\FG)$, if and only if $\FG\bot_\Lambda N_{\omega}$, i.e.
\beq\forall\,v\in N_{\omega}\quad\skp{\Lambda\FG}{v}_{\Lzqqpeom{}}=0.\mylabel{fgsenkrecht}\eeq
The solution can be chosen, such that $\EH\bot_\Lambda N_{\omega}$, i.e.
\beq\forall\,v\in N_{\omega}\quad\skp{\Lambda\EH}{v}_{\Lzqqpeom{}}=0.\mylabel{ehsenkrecht}\eeq
By this condition $\EH$ is uniquely determined
and the solution operator $\loesom\FG:=u$ is continuous.
More precisely, $\loesom$ maps $\Lzqqpesom\cap N_{\omega}^{\bot_\Lambda}$ to
$\big(\ronqtom\times\dqpetom\big)\cap N_{\omega}^{\bot_\Lambda}$
continuously for all $s,-t>1/2$. 
Here we denote the orthogonality corresponding to the
${}_\Lambda\Lzqqpeom{}$-scalar product by $\bot_\Lambda$.
\end{itemize}
\end{theo}

\begin{proof} 
We follow the proof of \cite[Theorem 2.10]{xmas}.
To show {\bf(i)}, i.e. the polynomial decay of any eigen-solution $\EH$, 
we only have to prove
$$\EH\in\Lzqqpegmehom$$
because of Theorem \ref{polyabklmax}, Remark \ref{bempolyabklmax}, 
the equation $M\EH=-\ie\omega\Lambda\EH$ and the inclusions
\beq\ol{\rot\ronqom{}}\subset\ronqpenom{},\quad\ol{\pdiv\dqpeom{}}\subset\dqnom{}.\mylabel{imsubsetkern}\eeq
Using the second part of the radiation condition we obtain some $t>-1/2$, such that
$$\lim_{\beta\to\infty}\norm{RE+H}_{\Lzqpet(\Bc_{r_{0}}\cap B_{\beta})}<\infty$$
holds true. We calculate
\begin{align*}
&\qquad\norm{RE+H}_{\Lzqpet(\Bc_{r_{0}}\cap B_{\beta})}^2\\
&=\norm{RE}_{\Lzqpet(\Bc_{r_{0}}\cap B_{\beta})}^2
+\norm{H}_{\Lzqpet(\Bc_{r_{0}}\cap B_{\beta})}^2
+2\ret\skp{\Phi RE}{H}_{\lzqpe(\Bc_{r_{0}}\cap B_{\beta})}
\end{align*}
with $\Phi:=\rho^{2t}$. Lemma \ref{techstrahl}, the differential equation and the
symmetry of $\eps$, $\mu$ yield
\begin{align*}
\skp{\Phi RE}{H}_{\lzqpe(\Bc_{r_{0}}\cap B_{\beta})}
&=\skp{\Psi\rot E}{H}_{\lzqpe(\Omega\cap B_{\beta})}
+\skp{\Psi E}{\pdiv H}_{\lzq(\Omega\cap B_{\beta})}\\
&=-\ie\omega\ub{\skp{\Psi\mu H}{H}_{\lzqpe(\Omega\cap B_{\beta})}}_{\in\rz}
+\ie\omega\ub{\skp{\Psi E}{\eps E}_{\lzq(\Omega\cap B_{\beta})}}_{\in\rz}.
\end{align*}
Thus, $\ret\skp{\Phi RE}{H}_{\lzqpe(\Bc_{r_{0}}\cap B_{\beta})}=0$ 
and hence, by means of the monotone convergence theorem
$H\in\Lzqpetom$ follows for $\beta\to\infty$.
Finally we get $E\in\Lzqgmehom$ using the first part of the radiation condition.

If {\bf(ii)} or {\bf(iii)} would be wrong
then there would exist a sequence of eigen-values
$(\omell)_{\ell\in\nz}\subset\rzon$ tending to $\omega$ 
and a sequence of eigen-forms
$(\EHell)_{\ell\in\nz}\subset N_{\omell}$,
such that $(\EHell)_{\ell\in\nz}$ is an ortho-normal system
with respect to the ${}_\Lambda\Lzqqpeom{}$-scalar product.
As an ortho-normal system $(\EHell)_{\ell\in\nz}$ converges in $\Lzqqpeom{}$ weakly to zero.
Moreover, by the differential equation $(\EHell)_{\ell\in\nz}$ is bounded in
$$\big(\ronqom{}\cap\eps^\me\dqnom{}\big)\times\big(\dqpeom{}\cap\mu^\me\ronqpenom{}\big).$$
Hence, from the {\sf MLCP} we can extract a subsequence $(\EHpiell)_{\ell\in\nz}$,
where $\pi:\nz\to\nz$ is strictly monotone,
converging for all $t<0$ in $\Lzqqpetom$ to $0$.
The latter is due to the weak convergence.
For $1\leq s\in\rz\ohne\pI$ Theorem \ref{polyabklmax} yields uniformly in $(\EHpiell)_{\ell\in\nz}$ and
$(\omepill)_{\ell\in\nz}$ the estimate
\begin{align*}
1=\skp{\Lambda\EHpiell}{\EHpiell}_{\Lzqqpeom{}}
&\leq c\norm{\EHpiell}_{\Lzqqpeom{}}^2
\leq c\norm{\EHpiell}_{\Lzqqpeom{s-1}}^2\\
&\leq c\norm{\EHpiell}_{\Lzqqpe{}(\Omega\cap B_\theta)}^2\xrightarrow{\ell\to\infty}0,
\end{align*}
which is a contradiction.

We prove {\bf(iv)}:
First of all \eqref{fgsenkrecht} is necessary
since we get for all eigen-forms $v\in N_{\omega}$ 
by their polynomial decay and Remark \ref{schaufelnbem}
$$\skp{\Lambda\FG}{v}_{\Lzqqpeom{}}
=\ie\skpb{(M+\ie\omega\Lambda)\EH}{v}_{\Lzqqpeom{}}
=-\ie\skpb{\EH}{\ub{(M+\ie\omega\Lambda)v}_{=0}}_{\Lzqqpeom{}}=0.$$
To show existence we now use Eidus' principle of limiting absorption. 
For that purpose let
$\FG\in\Lzqqpegehom$ with \eqref{fgsenkrecht}.
Moreover, let $(\sigma_\ell)_{\ell\in\nz}$ be a positive sequence tending
to zero and $(\FGell)_{\ell\in\nz}\subset\Lzqqpesom$ with some $s>1/2$
be a sequence satisfying $\FGell\bot_{\Lambda}N_{\omega}$,
such that $(\FGell)_{\ell\in\nz}$ converges to $\FG$ in $\Lzqqpesom$ as $\ell$ tends to infinity.
Defining non-real frequencies in the upper half plane $\omell\in\cz_+\ohne\rz$
with $\omell^2=\omega^2+\ie\sigma_\ell\omega$ and $\omell\to\omega$
we obtain $\lz$-solutions
$$\EHell:=\loes_{\omell}\FGell\in D(\calM)=\ronqom{}\times\dqpeom{}$$
solving the Maxwell problem $\Max(\Lambda,\omell,\FGell)$, i.e.
\beq(M_{\Lambda}-\omell)\EHell=\FGell.\mylabel{gleichunggrenzab}\eeq
Applying the orthogonal Hodge-Helmholtz decomposition we decompose
$\EHell$ and $\FGell$ according to Lemma \ref{hodgehelmholtzdeco}
$$\EHell=\EHellN+\EHellI,\FGell=\FGellN+\FGellI\in N(\calM)\oplus_{\Lambda}\ol{I(\calM)}.$$
Therefore,
$$(M_{\Lambda}-\omell)\EHellI-\FGellI=\omell\EHellN+\FGellN\in N(\calM)\cap\ol{I(\calM)}=\{0\}.$$
As orthogonal projections the sequences of forms $(\FGellN)_{\ell\in\nz}$, $(\FGellI)_{\ell\in\nz}$ 
converge in $\Lzqqpeom{}$. Hence, so does $(\EHellN)_{\ell\in\nz}$.
Let us assume the boundedness of $(\EHellI)_{\ell\in\nz}$ or 
\beq\forall\,t<-1/2\,\exists\,c>0\,\forall\,\ell\in\nz\quad
\norm{\EHell}_{\Lzqqpetom}\leq c\mylabel{annahme}\eeq
for a moment. At the end of the proof we will show by contradiction 
that in fact \eqref{annahme} holds.
Let $t'$ be such a $t$ with \eqref{annahme}. 
Then, $(\EHellI)_{\ell\in\nz}$ is bounded in $\Lzqqpeom{t'}$ and by the differential equation, i.e.
$$(M_{\Lambda}-\omell)\EHellI=\FGellI,$$
and by \eqref{imsubsetkern} even in 
$\big(\ronqom{t'}\cap\eps^\me\dqnom{t'}\big)\times\big(\dqpeom{t'}\cap\mu^\me\ronqpenom{t'}\big)$.
Hence, the {\sf MLCP} yields for an arbitrary $\tilde{t}<t'$ a subsequence 
$(\EHpiellI)_{\ell\in\nz}$
converging in $\Lzqqpeom{\tilde{t}}$ and even in $\ronqom{\tilde{t}}\times\dqpeom{\tilde{t}}$
by the differential equation. 
Therefore, the entire sequence $(\EHpiell)_{\ell\in\nz}$ converges in
$\ronqom{\tilde{t}}\times\dqpeom{\tilde{t}}$ to, let us say,
$$\EH\in\ronqom{\tilde{t}}\times\dqpeom{\tilde{t}},$$
which solves
$$(M_{\Lambda}-\omega)\EH=\FG.$$
With the polynomial decay of eigen-solutions and Remark \ref{schaufelnbem} we compute
for all eigen-forms $v\in N_{\omega}$ and all $\ell\in\nz$
$$0=\skp{\Lambda\FGpiell}{v}_{\Lzqqpeom{}}
=-\ie\skpb{\EHpiell}{(M+\ie\omega_{\pi\ell}\Lambda)v}_{\Lzqqpeom{}}
=\ub{(\ol{\omega}_{\pi\ell}-\omega)}_{\neq0}\skp{\Lambda\EHpiell}{v}_{\Lzqqpeom{}}.$$
Consequently, $\skp{\Lambda\EHpiell}{v}_{\Lzqqpeom{}}=0$.
Since $\skp{\,\cdot\,}{\Lambda v}_{\Lzqqpeom{}}$ is a continuous linear
functional on $\Lzqqpeom{\tilde{t}}$ for all $v\in N_{\omega}$ we obtain
\beq\forall\,v\in N_{\omega}\quad\skp{\Lambda\EH}{v}_{\Lzqqpeom{}}=0.\mylabel{EHsenkrecht}\eeq
Now, we pick some $t<-1/2$. Then, we get by Lemma \ref{apriorimaxwell} 
constants $\hat{t}>-1/2$ and $c,\theta>0$, such that by the monotone convergence theorem
$$\norm{\EH}_{\rqtom\times\dqpetom}
+\normb{(S+1)\EH}_{\Lzqqpeom{\hat{t}}}
\leq c\big(\norm{\FG}_{\Lzqqpesom}
+\norm{\EH}_{\Lzqqpe{}(\Omega\cap B_\theta)}\big)$$
holds. 
Therefore, $\EH\in\ronqkmehom\times\dqpekmehom$ 
and $\EH$ satisfies the radiation condition, i.e. 
$(S+1)\EH\in\Lzqqpegmehom$.
In other words, $\EH$ solves $\Max(\Lambda,\omega,\FG)$.

We note that this proves the principle of limiting absorption to hold.
The choice $\FGell:=\FG$ for all $\ell\in\nz$ yields the existence of a solution of
$\Max(\Lambda,\omega,\FG)$ and this one is unique because of \eqref{EHsenkrecht}.

Moreover, for $-t,s>1/2$ the solution operator $\loesom$ maps
$D_s(\loesom)$ to $I_t(\loesom)$, where 
$$D_s(\loesom):=\Lzqqpesom\cap N_{\omega}^{\bot_\Lambda},\quad
I_t(\loesom):=\big(\ronqtom\times\dqpetom\big)\cap N_{\omega}^{\bot_\Lambda},$$
continuously. This follows by the closed graph theorem 
because $D_s(\loesom)$ and $I_t(\loesom)$ are Hilbert spaces 
by the polynomial decay of eigen-solutions and $\loesom$ is closed.
The latter assertion is a consequence of Lemma \ref{apriorimaxwell} 
and the monotone convergence theorem.

Finally, it remains to contradict the contrary assumption to \eqref{annahme}. 
To this end, let $t<-1/2$ and
$(\EHell)_{\ell\in\nz}\subset\ronqtom\times\dqpetom$ be a sequence with
$\norm{\EHell}_{\Lzqqpetom}\to\infty$.
Defining the normalized forms
$$\EHsell:=\norm{\EHell}_{\Lzqqpetom}^\me\EHell,\quad
\FGsell:=\norm{\EHell}_{\Lzqqpetom}^\me\FGell$$
we have $\norm{\EHsell}_{\Lzqqpetom}=1$ for all $\ell\in\nz$
and $\norm{\FGsell}_{\Lzqqpesom}\to0$. Moreover, the equation
$$(M_{\Lambda}-\omell)\EHsell=\FGsell$$
holds. Following the arguments above we obtain a subsequence $(\EHspiell)_{\ell\in\nz}$
converging in $\Lzqqpeom{\tilde{t}}$ with $\tilde{t}<t$ towards
$\EHs\in N_{\omega}^{\bot_\Lambda}$, which solves $\Max(\Lambda,\omega,0)$.
Hence, $\EHs=0$ and Lemma \ref{apriorimaxwell} yields constants $c,\theta>0$
independent of $\sigma_{\pi\ell}$, $\FGspiell$ or $\EHspiell$, such that
$$1=\norm{\EHspiell}_{\Lzqqpetom}
\leq c\big(\ub{\norm{\FGspiell}_{\Lzqqpesom}}_{\to0}
+\ub{\norm{\EHspiell}_{\Lzqqpe{}(\Omega\cap B_\theta)}}_{\to0}\big)$$
holds true; a contradiction.
\end{proof}

The polynomial decay of eigen-solutions proved above and Theorem \ref{expabkl} yield:

\begin{cor}\mylabel{expdecayeigen}
Let $\omega\in\rzon$ and $\Lambda$ be $\tau$-$\pc{2}$-admissible
with $\tau>1$. Then, any eigen-solution $\EH\in N_{\omega}$
decays exponentially, i.e.
$$\e^{tr}\EH\in\big(\ronqom{}\cap\eps^\me\dqom{}\big)\times\big(\dqpeom{}\cap\mu^\me\ronqpeom{}\big)
\cap\qh{2}{q,q+1}{}{}(\Bc_{r_{0}+1})$$
holds for all $t\in\rz$.
\end{cor}

\begin{rem}\mylabel{polynomialdecayrem}
The polynomial resp. exponential decay of eigen-solutions holds for arbitrary exterior domains $\Omega$,
i.e. $\Omega$ does not need to have the {\sf MLCP}.
\end{rem}

\begin{rem}\mylabel{bemfredholmRellich}
If the medium is homogeneous and isotropic in the exterior of some ball, i.e.
$\supp\Lambdad\cup(\rN\ohne\Omega)\subset B_{\theta}$
for some $\theta>0$, then
$$\EH=0\quad\text{in }\Bc_\theta\,$$
for all $\omega\in\rzon$ and $\EH\in N_{\omega}$. 
Because in this case $\EH$ solves componentwise Helmholtz' equation $(\Delta+\omega^2)\EH=0$
in $\Bc_{\theta}$ and therefore, by Rellich's estimate \cite{rellich} or \cite[p. 59]{leisbuch},
must vanish in $\Bc_{\theta}$. 
If the principle of unique continuation holds for our Maxwell system under consideration
then even
$$N_{\omega}=\{0\}.$$
\end{rem}

Using the a priori estimate of the limiting absorption principle and some indirect arguments
followed by the 'trivial' decomposition of $\Lzqsom$ from \cite[Lemma 4.1]{paulydeco}
we are able to prove stronger estimates for the solution operator $\loesom$ 
as the ones given in Theorem \ref{fredholm} (iv).
We only note the results here.
For this, let $\om$ possess the {\sf MLCP} and let $s,-t>1/2$ as well as 
$K$ be a compact subset of $\czp\ohne\{0\}$.
Furthermore, let $\Lambda$ be $\tau$-admissible with some $\tau>1$.

\begin{lem}\mylabel{loesomabscheins}
There exist constants $c,\theta>0$ and some $\hat{t}>-1/2$, such that
$$\norm{\loesom\FG}_{\rqtom\times\dqpetom}
+\normb{(S+1)\loesom\FG}_{\Lzqqpeom{\hat{t}}}
\leq c\big(\norm{\FG}_{\Lzqqpesom}
+\norm{\loesom\FG}_{\Lzqqpe{}(\Omega\cap B_{\theta})}\big)$$
holds for all $\omega\in K$ and $\FG\in\Lzqqpesom\cap N_{\omega}^{\bot_{\Lambda}}$.
\end{lem}

\begin{cor}\mylabel{loesomabschzwei}
Let $K\cap P=\emptyset$. Then there exist constants $c>0$ and $\hat{t}>-1/2$, 
such that the estimate
$$\norm{\loesom\FG}_{\rqtom\times\dqpetom}
+\normb{(S+1)\loesom\FG}_{\Lzqqpeom{\hat{t}}}
\leq c\norm{\FG}_{\Lzqqpesom}$$
holds true for all $\omega\in K$ and $\FG\in\Lzqqpesom$. 
In particular, the solution operator
$\loesom$ mapping $\Lzqqpesom$ to $\ronqtom\times\dqpetom$
is equi-continuous with respect to $\omega\in K$.
\end{cor}

\begin{theo}\mylabel{stetigefortsetzungoperatornorm}
Let $K\cap P=\emptyset$. Then, the mapping
$$\Abb{\loes}{K}{\mathsf{B}\big(\Lzqqpesom,\ronqtom\times\dqpetom\big)}{\omega}{\loesom}$$
is uniformly continuous. 
Here we denote the set of bounded linear operators from some normed space $X$
to some normed space $Y$ by $\mathsf{B}(X,Y)$.
\end{theo}

\appendix

\section{Appendix: classical equations}

We want to point out briefly, which classical equations are covered by our generalized approach.
Since the relation between the differential form calculus and the classical vector calculus 
is very well known we directly translate our equations \eqref{shortproblem}, 
i.e. in the longer version \eqref{Maxglrotdiv}-\eqref{MaxRadCon},
into terms of vector analysis.

\begin{description}
\item $q=0$: $E$, $F$ are scalar functions and $H$, $G$ are vector fields. 
We have the equations of linear acoustics with Dirichlet boundary condition in first order form
\begin{align*}
\text{div}H+\ie\omega\eps E&=-\ie\eps F,\quad\text{grad}E+\ie\omega\mu H=-\ie\mu G&&\text{in }\om,\\
E&=0&&\text{on }\dom
\end{align*}
and $\xi\cdot H+E$, $E\xi+H$ decay at infinity, where $\xi(x):=x/r(x)$.
Using the differential equation for $H$, i.e. 
$$-\ie\omega(\xi\cdot H+E)=\ub{\xi\cdot\text{grad}E}_{=\p_{r}E}-\ie\omega E+\dots,$$
we see that we get Sommerfeld's classical radiation condition, i.e.
$(\p_{r}-\ie\omega)E$ decays at infinity.
\item $q=N-1$: $E$, $F$ are vector fields and $H$, $G$ are scalar functions. 
We have the equations of linear acoustics with Neumann boundary condition in first order form
\begin{align*}
\text{grad}H+\ie\omega\eps E&=-\ie\eps F,\quad\text{div}E+\ie\omega\mu H=-\ie\mu G&&\text{in }\om,\\
\nu\cdot E&=0&&\text{on }\dom
\end{align*}
and $H\xi+E$, $\xi\cdot E+H$ decay at infinity. 
Here $\nu$ is the unit normal vector pointing outwards. 
Now, using the differential equation for $E$ we get again Sommerfeld's radiation condition.
Note that by the differential equation the homogeneous boundary condition is equivalent to
the inhomogeneous Neumann boundary condition
$$\p_{\nu}H=-\ie\omega\nu\cdot\epsd(E-f)-\ie\nu\cdot F\quad\text{on }\dom.$$
\item $q=N$: $E$, $F$ are scalar functions and $H$, $G$ vanish. 
We have trivial equations
$$\omega E=-F\text{ in }\om,\quad
E\text{ decays at infinity}.$$
\item $q=1$, $N=3$: $E$, $F$, $H$, $G$ are all vector fields. 
We have the classical Maxwell equations with homogeneous electric boundary condition,
i.e. $\dom$ is a perfect conductor,
\begin{align*}
-\text{curl}H+\ie\omega\eps E&=-\ie\eps F,\quad\text{curl}E+\ie\omega\mu H=-\ie\mu G&&\text{in }\om,\\
\nu\times E&=0&&\text{on }\dom
\end{align*}
and $-\xi\times H+E$, $\xi\times E+H$ decay at infinity.
The latter are the classical Silver-M\"uller radiation conditions for Maxwell's equations.
\end{description}

\section{Appendix: some calculations}

Let $\tr$ and $\sw$ denote the formal trace and the formal swapping operator, respectively, i.e. 
$$\tr\zmat{A}{a}{b}{B}=A+B,\quad \sw\zmat{A}{a}{b}{B}=\zmat{B}{a}{b}{A}.$$
We note
$$S^2=\zmat{TR}{0}{0}{RT},\quad 1-S^2=\zmat{RT}{0}{0}{TR}$$
and formally $1-S^2=\sw S^2$ and $\tr S^2=\tr\sw S^2=RT+TR=1$.
Moreover,
$$SM=\zmat{T\rot}{0}{0}{R\pdiv},\quad MS=\zmat{\pdiv R}{0}{0}{\rot T}$$
and formally $\tr SM=T\rot+R\pdiv$, $\tr MS=\rot T+\pdiv R$.
By a straight forward computation we obtain:

\begin{lem}\mylabel{gewichttausch}
Let $m\in\rz$.
\begin{align*}
\text{\rm\bf(i)}&&C_{\p_{n},r^m}&=\frac{m}{r}r^{-1}\id_{n}r^{m},\quad 
C_{\Delta,r^m}=\frac{m}{r}(2\p_r-\frac{m+2-N}{r})r^m\\
\text{\rm\bf(ii)}&&C_{\rot,r^m}&=\frac{m}{r}Rr^m,\quad 
C_{\pdiv,r^m}=\frac{m}{r}Tr^m,\quad C_{M,r^m}=\frac{m}{r}Sr^m\\
\text{\rm\bf(iii)}&&C_{\pdiv\rot,r^m}
&=\frac{m}{r}(\pdiv R+T\rot-\frac{m+1}{r}TR)r^m,\quad
C_{\rot\pdiv,r^m}=\frac{m}{r}(\rot T+R\pdiv-\frac{m+1}{r}RT)r^m\\
\text{\rm\bf(iv)}&&C_{M^2,r^m}
&=\frac{m}{r}(MS+SM-\frac{m+1}{r}S^2)r^m\\
\text{\rm\bf(v)}&&C_{\Box,r^m}
&=\frac{m}{r}(\sw MS+\sw SM-\frac{m+1}{r}\sw S^2)r^m\\
\text{\rm\bf(vi)}&&C_{\Delta,r^m}
&=\frac{m}{r}(\tr MS+\tr SM-\frac{m+1}{r})r^m\\
\end{align*}
Looking at the two formulas for $C_{\Delta,r^m}$ 
we get the nice relation
$$\rot T+\pdiv R+T\rot+R\pdiv=\tr(SM+MS)=(N-1)/r+2\p_r.$$
Here $\p_r$ is meant componentwise in Cartesian coordinates.
\end{lem}

\begin{acknow}
The author is particularly indebted to his academic teachers Norbert Weck and Karl-Josef Witsch.
\end{acknow}

\end{document}